\documentclass[11pt,a4paper]{amsart}
\usepackage[T1]{fontenc}
\usepackage[utf8]{inputenc}
\usepackage{lmodern}
\usepackage{amsmath,amssymb,booktabs,needspace}
\usepackage{microtype,cite}
\DeclareMicrotypeSet{articletext}{encoding={T1},family={rm*,sf*}}
\UseMicrotypeSet[expansion]{articletext}
\usepackage[margin=27mm,headheight=12pt,headsep=12pt]{geometry}
\usepackage[hidelinks,pdfencoding=auto]{hyperref}
\hypersetup{
 pdftitle={A sphere-product counterexample to the additivity of homology nilpotency},
 pdfauthor={Paul-Eugène Parent},
 pdfsubject={Homology nilpotency, smooth manifolds, sphere products, and Koszul residues}}
\newtheorem{theorem}{Theorem}[section]
\newtheorem{proposition}[theorem]{Proposition}
\newtheorem{lemma}[theorem]{Lemma}
\newtheorem{corollary}[theorem]{Corollary}
\theoremstyle{definition}
\newtheorem{definition}[theorem]{Definition}
\theoremstyle{remark}

\numberwithin{equation}{section}
\newcommand{\Q}{\mathbb Q}
\newcommand{\M}{\mathcal M}
\newcommand{\Qf}{\mathcal Q}
\newcommand{\N}{\mathcal N}
\newcommand{\A}{\mathcal A}
\newcommand{\E}{\mathcal E}
\newcommand{\F}{\mathcal F}
\newcommand{\K}{\mathcal K}
\newcommand{\W}{\mathcal W}
\newcommand{\Pmod}{\mathcal P}
\newcommand{\m}{\mathfrak m}
\newcommand{\q}{\mathfrak q}
\newcommand{\res}{\mathfrak r}
\newcommand{\ellid}{\mathfrak l}
\newcommand{\hn}{\operatorname{Hnil}}
\newcommand{\cat}{\operatorname{cat}}
\newcommand{\Cat}{\operatorname{Cat}}
\newcommand{\nilh}{\operatorname{nil}_{h}}
\newcommand{\Ann}{\operatorname{Ann}}
\newcommand{\Span}{\operatorname{span}}
\newcommand{\ah}{\widehat a}
\newcommand{\ctr}{\mathbin{\lrcorner}}
\newcommand{\eps}{\varepsilon}
\newcommand{\doi}[1]{\href{https://doi.org/#1}{\nolinkurl{doi:#1}}}
\allowdisplaybreaks[1]
\AtBeginDocument{\flushbottom}

\title[Homology nilpotency and sphere products]{A sphere-product counterexample to the additivity of homology nilpotency}
\author[Paul-Eugène Parent]{Paul-Eugène Parent$^{*}$}
\address{Department of Mathematics and Statistics, University of Ottawa,
Ottawa, Ontario K1N 6N5, Canada}
\email{\href{mailto:pparent@uottawa.ca}{pparent@uottawa.ca}}
\date{November 1st,2026}
\thanks{$^{*}$This research was supported by the Natural Sciences and Engineering
Research Council of Canada (NSERC), grant no.~149521.}

\begin{document}
\begin{abstract}
Homology nilpotency need not increase under a product with a sphere.
We construct a simply connected, closed smooth manifold $M$ of dimension
$53$, rationally elliptic with $\dim_{\Q}\bigl(\pi_*(M)\otimes\Q\bigr)=7$, such that
$\cat_0(M)=10$ and $\hn(M\times S^2)=\hn(M)=11$.
The inequality $\hn(X\times Y)\leq\hn(X)+\hn(Y)$ can therefore be
strict even for a product of a closed manifold with $S^2$ as the second
factor.
The upper bound follows from a recognition theorem for acyclic ideals in
minimal Sullivan algebras of the form $\Lambda(a)\otimes\K$, with
differential $\partial+a\theta$, whose pure factor $\K$ has Poincaré
duality in cohomology with orientation in Koszul degree two.
Polynomial residues of elements in the radical of an orientation pairing
control ambient boundaries and the exactness of contraction by $a$.
Poincaré duality and a comparison of two short exact sequences then
establish acyclicity. For the example, five corrected elements verify the
criterion for a two-parameter family of orientation functionals.
The quadratic sphere relation accounts for a nonzero second differential
in a sphere-weight filtration. Barge--Sullivan realization preserves the
model and allows all rational Pontryagin classes to vanish.
\end{abstract}
\maketitle
\tableofcontents

\section{Introduction}\label{sec:intro}
Product formulas distinguish closely related invariants in rational
homotopy theory. For Lusternik--Schnirelmann category, the product question
is classical. We use the normalized convention $\cat(*)=0$ and
$\cat(S^k)=1$, and write $\cat_0(X)=\cat(X_{\Q})$ for rational category.
Ganea asked whether
\begin{equation}\label{eq:ganea}
 \cat(X\times S^k)=\cat(X)+1
\end{equation}
holds for finite complexes $X$ and $k\geq1$
\mbox{\cite[Problem~2]{Ganea}}. Iwase's counterexamples show that it does not
hold in general \mbox{\cite{Iwase}}. The surrounding theory is developed in
\mbox{\cite{CLOT}}.

The rational counterpart of \eqref{eq:ganea} does hold. Halperin and
Lemaire introduced module category $\operatorname{Mcat}$ by replacing the algebra retraction
condition for rational category with a retraction of differential modules
\mbox{\cite[p.~138]{HL}}. Jessup proved the sphere-product formula for this
invariant \mbox{\cite[Corollary~3.3, p.~62]{Jessup}}, and Hess proved
$\operatorname{Mcat}=\cat_0$ \mbox{\cite[\mbox{Theorem}~0]{Hess}}. Thus
$\cat_0(X\times S^k)=\cat_0(X)+1$ for simply connected spaces of finite
rational type and $k\geq2$. Félix, Halperin, and Lemaire proved the stronger
product theorem \mbox{\cite[\mbox{Theorem}~1]{FHL}}:
\begin{equation}\label{eq:cat-product}
 \cat_0(X\times Y)=\cat_0(X)+\cat_0(Y).
\end{equation}

This paper concerns the additional constraint of realizing a rational
homotopy type by a quotient of its minimal model. We abbreviate
\emph{commutative differential graded algebra} to cdga; differentials have
cohomological degree $+1$. The nilpotency length of an augmented cdga is
the least $n$ for which its augmentation ideal has zero $(n+1)$st power.
For a simply connected minimal Sullivan algebra $\M=(\Lambda V,D)$, put
$\m_{\M}=\Lambda^{\geq1}V$. The \emph{homology nilpotency} of $\M$ is
\begin{equation}\label{eq:hn-definition}
 \hn(\M)=\inf\{n\in\mathbb Z_{\geq0}:\m_{\M}^{n+1}\subset J,
       \ J\text{ an acyclic differential ideal of }\M\},
\end{equation}
where $\inf\varnothing=+\infty$. Ideals are graded by cohomological
degree; no auxiliary homogeneity is required. Equivalently, $\hn(\M)$
is the least nilpotency length of a quotient whose projection from $\M$
is a quasi-isomorphism. For a space, we use its minimal rational model.

Carrasquel-Vera introduced homology nilpotency for an ideal in a cdga
\mbox{\cite[Definition~20]{Carrasquel}}; \eqref{eq:hn-definition} is
$\hn(\m_{\M})$ in that notation. The invariant provides an upper bound
for rational sectional category: for a surjective cdga model, the
homology nilpotency of the kernel bounds the associated sectional category
\mbox{\cite[Proposition~21]{Carrasquel}}. Applied to the augmentation ideal, this
leads to a product problem that retains the distinction between a quotient
of the minimal model and an arbitrary model of the same rational type.

Let $\nilh(\M)$ denote the least nilpotency length of any cdga
quasi-isomorphic to $\M$, without requiring a quotient map from $\M$.
The subscript $h$ stands for \emph{homotopical}; we call $\nilh$ the
homotopical nil-length. Write $\Cat_0(X)=\Cat(X_{\Q})$ for rational
strong category, where $\Cat(X)$ is the least $n$ such that some space of
the homotopy type of $X$ has a cover by $n+1$ contractible open sets.
Cornea identifies rational cone length with homotopical nil-length
\mbox{\cite{CorneaCone}} and with rational strong category
\mbox{\cite[Definition~2, \mbox{Theorem}~1.1 and Corollary~2, p.~837]{Cornea}}.
With $e_0$ denoting the rational Toomer invariant, the standard comparisons
are
\begin{equation}\label{eq:comparisons}
 e_0(X)\leq\cat_0(X)\leq\Cat_0(X)=\nilh(X)\leq\hn(X).
\end{equation}
The first four terms coincide for simply connected rational Poincaré
duality spaces by \mbox{\cite[\mbox{Theorem}~1.1, p.~744]{CFL}} and
\mbox{\cite[\mbox{Theorem}~3]{FHL}}.
We call $\hn(X)-\nilh(X)$ the \emph{quotient defect} when these invariants
are finite. On Poincaré duality spaces this is $\hn(X)-\cat_0(X)$.

For a formal minimal algebra the quotient defect is zero. Formality gives
a quasi-isomorphism $\M\to(H(\M),0)$ \mbox{\cite{FHT}}. This map is surjective:
every element of $H(\M)$ is the image of a cocycle, because the induced
map in cohomology is onto. The acyclic kernel gives an upper bound for
$\hn(\M)$ by the cup length of $H(\M)$. Together with the cup-length
lower bound for category and \eqref{eq:comparisons}, this proves that the
cup length, $\cat_0(\M)$, and $\hn(\M)$ coincide. Nonformality is
therefore necessary for the phenomenon studied here.

Tensoring quotient models gives
\begin{equation}\label{eq:hn-product}
 \hn(X\times Y)\leq\hn(X)+\hn(Y)
\end{equation}
(Proposition~\ref{prop:quotient}). Our example shows that one sphere factor
can leave homology nilpotency unchanged, even for closed smooth manifolds.

\begin{theorem}\label{thm:main}
There exists an oriented, simply connected, closed smooth manifold $M$ of
dimension $53$, rationally elliptic with $\dim_{\Q}\bigl(\pi_*(M)\otimes\Q\bigr)=7$, such that
\[
 \cat_0(M)=10,\qquad
 \hn(M)=\hn(M\times S^2)=\cat_0(M\times S^2)=11.
\]
The manifold can be chosen with $p_i(TM)=0\in H^{4i}(M;\Q)$ for every
$i\geq1$. In particular,
\begin{equation}\label{eq:strict}
 \hn(M\times S^2)=11<12=\hn(M)+\hn(S^2).
\end{equation}
\end{theorem}

The quotient defect of $M$ is one, and a product with $S^2$ removes it.
A related operation is stabilization by a degreewise finite wedge of
rational spheres, which can lower homology nilpotency to rational category
\mbox{\cite[\mbox{Theorem}~1.2, p.~8]{Parent}}. Here the operation is a product with
one sphere: homology nilpotency remains fixed and category increases.
These two operations have different effects on the rational homotopy type.

\subsection{The mechanism of the proof}
The minimal model of the manifold in Theorem~\ref{thm:main} is
$\N=(\Lambda(x,z,a,b_0,b_1,b_2,b_3),D)$, with degrees
$(2,6,3,11,13,15,17)$ and differential
\begin{equation}\label{eq:model-intro}
 \begin{gathered}
 Dx=Da=0,\qquad Dz=ax^2,\qquad Db_0=x^6,\\
 Db_1=x^4z+ab_0,\qquad Db_2=x^2z^2+2ab_1,\qquad
 Db_3=z^3+3ab_2.
 \end{gathered}
\end{equation}
We denote its rational realization by $N$ and its smooth realization by
$M$. The model of the product with $S^2$ is
$\A=\N\otimes\Lambda(u,v)$, where $|u|=2$, $|v|=3$, $Du=0$, and
$Dv=u^2$. Rational surgery realizes $\N$ in its formal dimension; the
nilpotency calculation uses these degrees throughout.

The central step is the construction of an acyclic ideal containing
$\m_{\A}^{12}$. The algebra has the form
$\A=\Lambda(a)\otimes\K$, with $D=\partial+a\theta$, where $\partial$
is a pure Koszul differential and $\partial\theta=\theta\partial$.
We choose a closed orientation functional $\eps$ and construct
\[
 \E=\F(\W_{10})+\A\ctr\eps\subset\A^\vee,
 \qquad J=\E^\perp.
\]
Here $\A\ctr\eps$ denotes the cyclic contraction module generated by
$\eps$ in the graded dual $\A^\vee$:
\[
 \A\ctr\eps=\Span_{\Q}\{r\ctr\eps:r\in\A\},
 \qquad (r\ctr\eps)(s)=\eps(sr)\quad(r,s\in\A).
\]
Thus contraction is defined by multiplying the argument on the right
before applying $\eps$.
The space $\W_{10}$ consists of polynomial coefficient functionals supported
on monomials of length at most ten, and $\F(\W_{10})$ is its differential contraction
module after adjoining the dual copy of $a$. Every functional in $\E$
has support of length at most eleven, so $\m_{\A}^{12}\subset J$.

The recognition theorem of Section~\ref{sec:recognition} proves that
$\E\to\A^\vee$ is a quasi-isomorphism. Its key input is a polynomial
ideal. For an element $\kappa+a\eta$ in the radical $O_\eps$ of the
orientation pairing, with $\kappa$ of Koszul degree one and $\eta$ of
Koszul degree two, retain the polynomial $\partial\kappa$. These
\emph{residues} form an ideal $\res_\eps\subset\Q[x,z,u]$.
The inclusion
\[
 (x,z,u)^{11}\subset\res_\eps
\]
forces every ambient dual boundary lying in $V_1+I_\Psi\ctr\Phi$ to
have a scalar primitive in $\W_{10}$; these spaces are defined in
Section~\ref{sec:recognition-proof}. This saturation property, together
with pure Poincaré duality, makes contraction by $a$ exact on $\E$.
A comparison of two short exact sequences, whose outside terms model
$\K^\vee$, gives the required quasi-isomorphism.

The residue condition is verified by five elements organized around the
adjacent syzygies of the four cubics $t^3,t^2z,tz^2,z^3$, with $t=x^2$.
The relation $Dv=u^2$ corrects the first primitive, while the last Koszul
syzygy provides a second correction. Their compatibility with every
multiplication is supplied by the residue theorem. The family of
orientation functionals has two parameters, and the construction works
exactly when both are nonzero. A sphere-weight filtration also detects a
specific obstruction that survives through weight one and dies at weight
two.

Section~\ref{sec:prelim} develops the quotient and duality tools.
Section~\ref{sec:model} determines the invariants of the factor, and
Section~\ref{sec:smooth} establishes smooth realization.
Section~\ref{sec:product} introduces the product and its orientation
functionals. Sections~\ref{sec:recognition} and~\ref{sec:recognition-proof}
state and prove the residue criterion; Section~\ref{sec:application}
applies it to the product. Section~\ref{sec:mechanism} treats the two
parameters and the filtered obstruction, and Section~\ref{sec:consequences}
gives repeated products and further questions.

\section{Quotient models and duality}\label{sec:prelim}
We work over $\Q$. All cdgas are connected, augmented, and concentrated
in nonnegative cohomological degrees. The symbol $\Lambda$ denotes the
free graded-commutative algebra, polynomial on even generators and
exterior on odd generators. Our minimal Sullivan algebras are simply
connected and finite-dimensional in each degree. We use Sullivan models
and their realization in the sense of \mbox{\cite{FHT,Sullivan}}.
A simply connected space is \emph{rationally elliptic} when both its
rational cohomology and its rational homotopy
$\pi_*(X)\otimes\Q=\bigoplus_{j\geq2}\pi_j(X)\otimes\Q$ are
finite-dimensional. For a minimal model
$(\Lambda V,D)$ this means that $V$ and $H(\Lambda V,D)$ are
finite-dimensional. All category values of manifolds in this paper are
rational category values.

For $\M=(\Lambda V,D)$, one has
$\m_{\M}^n=\Lambda^{\geq n}V$. Word length counts generators, with
multiplicity, and is distinct from cohomological degree. The nilpotency
length of a connected augmented cdga $C$ is the least $n$ such that
$(C^+)^{n+1}=0$, or $+\infty$ when no such integer exists.

\begin{proposition}\label{prop:quotient}
For a minimal Sullivan algebra $\M$, one has $\hn(\M)\leq n$ if and only
if there is a surjective quasi-isomorphism $\M\twoheadrightarrow C$ with
$(C^+)^{n+1}=0$. Moreover,
\[
 \hn(\M\otimes\M')\leq\hn(\M)+\hn(\M')
\]
for simply connected minimal Sullivan algebras $\M$ and $\M'$.
\end{proposition}
\begin{proof}
The sequence $0\to J\to\M\to\M/J\to0$ shows that the projection is a
quasi-isomorphism precisely when $J$ is acyclic. An acyclic ideal in a
connected cdga is proper, hence lies in the augmentation ideal. The
quotient has nilpotency length at most $n$ precisely when
$\m_{\M}^{n+1}\subset J$.

For the product assertion it suffices to consider finite values
$p=\hn(\M)$ and $q=\hn(\M')$. Choose quotient quasi-isomorphisms
$\M\twoheadrightarrow C$ and $\M'\twoheadrightarrow C'$ with
$(C^+)^{p+1}=0$ and $(C'^+)^{q+1}=0$. Their tensor product is a
surjective quasi-isomorphism by the Künneth theorem. A product of
$p+q+1$ elements of the augmentation ideal of $C\otimes C'$ expands into
terms containing either $p+1$ positive-degree factors from $C$ or $q+1$
from $C'$. Every term is zero. The first assertion applies to the product
minimal model.
\end{proof}

We shall eliminate odd generators whose differentials are regular. The
following formulation does not require a differential subalgebra on the
remaining generators. This generality is needed for the second
elimination in Proposition~\ref{prop:finite-quotient}.

\begin{lemma}[Regular elimination]\label{lem:regular}
Suppose the underlying graded algebra of a cdga is
$R_0\otimes\Lambda(w)$, with $w$ odd and $Dw=f\in R_0$. If
multiplication by $f$ is injective on $R_0$, then the differential ideal
$(w,f)$ is acyclic.
\end{lemma}
\begin{proof}
Write $Dn=D_0n+w\Delta_w n$ for $n\in R_0$. Since $Df=0$, the ideal
$(w,f)$ is differential. An element of the ideal has the form $fn+wn'$.
If it is closed, the component of its differential without $w$ is
$f(D_0n+n')$. Regularity gives $n'=-D_0n$, and therefore
\[
 fn+wn'=fn-wD_0n=D(wn),\qquad wn\in(w,f).\qedhere
\]
\end{proof}

For the minimal model $(\Lambda(u,v),Dv=u^2)$ of $S^2$, the ideal
$(v,u^2)$ is acyclic and the quotient $\Lambda(u)/(u^2)$ has nilpotency
length one. Nonzero positive-degree cohomology excludes length zero, so
\begin{equation}\label{eq:sphere}
 \hn(S^2)=1.
\end{equation}

\subsection{Coefficient functionals and contraction}
For a graded vector space $A$ finite-dimensional in each degree, let
\[
 A^\vee=\bigoplus_{d\geq0}\operatorname{Hom}_{\Q}(A^d,\Q).
\]
A functional on $A^d$ is indexed here by its \emph{support degree} $d$.
The transpose $D^*\xi=\xi\circ D$ lowers this index by one. This
convention has the same cycles and boundaries as the usual signed dual
convention. We write $H(A^\vee)$ for the resulting dual cohomology, with
support degrees understood. Right contraction is
\begin{equation}\label{eq:contraction}
 (r\ctr\xi)(s)=\xi(sr).
\end{equation}
For ordered monomials, the sign in \eqref{eq:contraction} is the sign
obtained by putting $sr$ in the fixed order of its odd letters. If
$\xi D=0$ and $\xi$ has support degree $d$, the Leibniz rule gives
\begin{equation}\label{eq:contraction-differential}
 D^*(r\ctr\xi)=(-1)^{d-|r|}(Dr)\ctr\xi
\end{equation}
for homogeneous $r$. To check the sign, evaluate on $s$ of degree
$d-|r|-1$ and use $\xi(D(sr))=0$.

A functional has \emph{support of length at most $n$} when it vanishes on
all monomials of length greater than $n$. Contraction does not increase
support length. In a minimal Sullivan algebra, each nonzero term of $D$
increases word length by at least one; consequently $D^*$ decreases
support length by at least one.

\begin{lemma}[Dual construction of an acyclic ideal]\label{lem:dual}
Let $A$ be a connected cdga finite-dimensional in each degree, and let
$U\subset A^\vee$ be a finite-dimensional graded subspace stable under
$D^*$ and all contractions. Then
\[
 U^\perp=\{s\in A:\xi(s)=0\text{ for every }\xi\in U\}
\]
is a differential ideal and $(A/U^\perp)^\vee\cong U$ as complexes.
If $U\hookrightarrow A^\vee$ is a quasi-isomorphism, then $U^\perp$ is
acyclic.
\end{lemma}
\begin{proof}
For $s\in U^\perp$, $r\in A$, and $\xi\in U$, stability gives
$\xi(sr)=(r\ctr\xi)(s)=0$ and
$\xi(Ds)=(D^*\xi)(s)=0$. Thus $U^\perp$ is a differential ideal.
The finite-dimensional double-annihilator theorem in each degree
identifies the quotient dual with $U$. Graded duality is exact in each
degree, so a quasi-isomorphism $U\to A^\vee$ dualizes to a
quasi-isomorphism $A\to A/U^\perp$. The exact sequence of its kernel
gives the last assertion.
\end{proof}

\subsection{Orientation pairings and their radicals}
A connected finite-dimensional graded algebra $H$ has \emph{Poincaré
duality of formal dimension $d$} if $H^d=\Q\omega$, $H^{>d}=0$, and
multiplication followed by evaluation on $\omega$ gives perfect pairings
$H^j\otimes H^{d-j}\to\Q$. A closed functional on a cdga that induces
this evaluation is called an orientation functional.

We use two consequences of this pairing. First, if a finite-dimensional
cdga has a closed top-degree functional whose multiplication pairing is
perfect on cochains, the induced pairing on cohomology is perfect. The
Leibniz rule makes the differential adjoint to itself up to degree signs.
Boundaries are contained in the orthogonal complement of cycles, and
adjointness and rank--nullity give equality. This proves the assertion in
each complementary pair of degrees.

Second, if $H(A)$ has Poincaré duality and $\eps$ is a closed orientation
functional, contraction by cycles induces the duality isomorphism
\begin{equation}\label{eq:PD-contraction}
 H(A)\longrightarrow H(A^\vee),\qquad
 [r]\longmapsto[r\ctr\eps],
\end{equation}
with the support-degree shift and signs from
\eqref{eq:contraction-differential}. On cohomology this map is precisely
the perfect multiplication pairing. This is the differential-module
duality used in \mbox{\cite[Lemma~3.2, p.~746]{CFL}}.

The \emph{ideal of orphans} of $\eps$ is
\begin{equation}\label{eq:orphans}
 O_\eps=\{r\in A:\eps(sr)=0\text{ for every }s\in A\}.
\end{equation}
This is the terminology of Lambrechts--Stanley
\mbox{\cite[Definition~3.1 and Proposition~3.2, p.~501]{LS}}.
The Leibniz rule proves that \eqref{eq:orphans} is a differential ideal.
A closed element of $O_\eps$ is a boundary in $A$: otherwise Poincaré
duality would give a complementary cocycle with nonzero pairing.
The primitive need not belong to $O_\eps$, so this observation does not
imply acyclicity of the ideal.

Lambrechts--Stanley construct a Poincaré duality cdga model for any
simply connected rational Poincaré duality type
\mbox{\cite[\mbox{Theorem}~1.1]{LS}}. That existence result permits a change of
cdga model. Homology nilpotency imposes a different requirement: the
projection must start from the fixed minimal algebra. Our construction
will keep that algebra and choose a finite submodule of its dual.

\subsection{Category from the length of an orientation}
The rational Toomer invariant originates in the Moore spectral sequence
\mbox{\cite{Toomer}}. We use its word-length description
\mbox{\cite[Remark~9.3(1), p.~29]{FH}}; see also \mbox{\cite{FHT}}:
\[
 e_0(\M)=\sup\{n:\text{some cocycle }\gamma\in\m_{\M}^{n}
                           \text{ has }[\gamma]\ne0\}.
\]
The word-length characterization of rational category
\mbox{\cite[\mbox{Theorem}~4.7]{FH}} gives $e_0(\M)\leq\cat_0(\M)$.
For simply connected rational Poincaré duality spaces the two invariants
are equal \mbox{\cite[\mbox{Theorem}~3]{FHL}}.

\begin{lemma}\label{lem:toomer-trace}
Suppose $H(\M)$ has Poincaré duality of formal dimension $d$.
If an orientation class has a representative $\zeta\in\m_{\M}^{n}$,
and if a functional $\lambda:\M^d\to\Q$ satisfies
\[
 \lambda D=0,\qquad \lambda(\zeta)=1,\qquad
 \lambda(\m_{\M}^{n+1}\cap\M^d)=0,
\]
then $e_0(\M)=\cat_0(\M)=n$.
\end{lemma}
\begin{proof}
The representative gives $e_0(\M)\geq n$. A nonzero class with a cocycle
in $\m_{\M}^{n+1}$ would, by duality, have a complementary cocycle whose
product represents a nonzero multiple of the orientation. That product
still belongs to $\m_{\M}^{n+1}$, contradicting the three properties of
$\lambda$. Hence $e_0(\M)\leq n$, and the category assertion follows
from the theorem just recalled.
\end{proof}

\section{The elliptic factor}\label{sec:model}
The model $\N$ in \eqref{eq:model-intro} has Sullivan order
$x,a,z,b_0,b_1,b_2,b_3$. In uniform notation,
\[
 Db_i=x^{6-2i}z^i+iab_{i-1}\quad(0\leq i\leq3),\qquad b_{-1}=0.
\]
For $1\leq i\leq3$, the two terms have differentials
\[
 D(x^{6-2i}z^i)=iax^{8-2i}z^{i-1},\qquad
 D(iab_{i-1})=-iax^{8-2i}z^{i-1}.
\]
They cancel, so $D^2=0$. All nonzero differentials are decomposable,
and all generators have degree at least two. Thus $\N$ is simply
connected and minimal. The fixed order of its odd letters is
$a,b_0,b_1,b_2,b_3$.

\subsection{A twisted Koszul complex and its syzygies}
Define the pure differential and an even derivation of degree $-2$ on
$\K_N=\Lambda(x,z,b_0,b_1,b_2,b_3)$ by
\[
 \begin{gathered}
 \partial x=\partial z=0,\qquad \partial b_i=x^{6-2i}z^i,\\
 \theta x=0,\qquad\theta z=x^2,\qquad\theta b_0=0,\qquad
 \theta b_i=ib_{i-1}\quad(1\leq i\leq3).
 \end{gathered}
\]
The derivation $\theta$ preserves the number of letters $b_i$, a number
that we call the Koszul degree. Since $\theta(\partial b_i)=\partial(\theta b_i)$,
\begin{equation}\label{eq:D-split}
 \partial\theta=\theta\partial,\qquad
 D=\partial+a\theta\text{ on }\K_N,\qquad Da=0.
\end{equation}
Putting $t=x^2$ identifies the pure part with the Koszul complex of the
four generators of $(t,z)^3$, extended from $\Q[t]$ to $\Q[x]$.
Its adjacent syzygies are
\begin{equation}\label{eq:syzygies}
 \begin{gathered}
 \sigma_i=x^2b_{i+1}-zb_i\quad(0\leq i\leq2),\qquad
 \partial\sigma_i=0,\\
 D\sigma_0=0,\qquad D\sigma_1=a\sigma_0,\qquad
 D\sigma_2=2a\sigma_1.
 \end{gathered}
\end{equation}
Their cohomological degrees are $17,19,21$.

The same structure describes the rational fibration with base $S^3_\Q$
modeled by
\[
 (\Lambda(a),0)\longrightarrow(\Lambda(a)\otimes\K_N,D)
                        \longrightarrow(\K_N,\partial).
\]
Successive adjunctions in the Sullivan order encode rational Postnikov
stages: the differential of a generator of degree $d$ represents the
corresponding $k$-invariant in degree $d+1$ \mbox{\cite[Sections~8 and~10]{Sullivan}}.
In these coordinates, the relation $\partial\theta=\theta\partial$
ensures that the successive representatives are closed. The terms with
$a$ in \eqref{eq:model-intro} are therefore linked by the single
derivation $\theta$. Section~\ref{sec:postnikov} returns to these
Postnikov data to distinguish the retained twist from the quotient
obstruction altered by the sphere factor.

\subsection{A finite quotient and its Poincaré pairing}
\begin{proposition}\label{prop:finite-quotient}
The ideal
\[
 I_{\mathrm{reg}}=(b_0,x^6,b_3,z^3+3ab_2)
\]
is acyclic. The quotient
\begin{equation*}
 \Qf=\N/I_{\mathrm{reg}}
   =\frac{\Lambda(x,z,a,b_1,b_2)}{(x^6,z^3+3ab_2)}
\end{equation*}
is finite-dimensional, has nilpotency length eleven, and has Poincaré
duality in cohomology of formal dimension $53$.
\end{proposition}
\begin{proof}
First eliminate $(b_3,z^3+3ab_2)$ using Lemma~\ref{lem:regular}.
The element $z^3+3ab_2$ is even, hence central, and is monic in $z$ over
$\Lambda(x,a,b_0,b_1,b_2)$. Multiplication by it is therefore injective. After this elimination the underlying algebra is free over $\Q[x]$,
so multiplication by $x^6$ is injective. Eliminate $(b_0,x^6)$ next.
The composite projection is a quasi-isomorphism, proving acyclicity of
its kernel.

A vector-space basis of $\Qf$ is
\begin{equation}\label{eq:Q-basis}
 x^p z^q a^e b_1^f b_2^g,
 \qquad 0\leq p\leq5,\quad 0\leq q\leq2,\quad
 e,f,g\in\{0,1\}.
\end{equation}
In particular $\Qf$ has dimension $144$. To compute the nilpotency length, consider
monomials before reducing the $z$-exponent. If that exponent is at most
two, a nonzero monomial has length at most $5+2+3=10$. If it is between
three and five, the relation $z^3=-3ab_2$ forbids another $a$ or $b_2$,
and the length is at most $5+5+1=11$. Exponents at least six give zero.
Thus $\m_{\N}^{12}\subset I_{\mathrm{reg}}$. The length-eleven word
$x^5z^5b_1$ is nonzero, since its image is $3ax^5z^2b_1b_2$.

For duality, the algebra
$\mathcal B=\Lambda(x,a,b_1,b_2)/(x^6)$ has the perfect pairing given by the
coefficient of $ax^5b_1b_2$. The extension
$\Qf=\mathcal B[z]/(z^3+3ab_2)$ is free over $\mathcal B$ with basis $1,z,z^2$.
Extraction of the $z^2$-coefficient gives a perfect relative pairing:
its matrix has ones on the antidiagonal and zeros elsewhere. The composite trace extracts the coefficient of
$ax^5z^2b_1b_2$, of degree $53$.

Every other basis monomial in \eqref{eq:Q-basis} has degree at most $51$,
so $\Qf^{52}=0$ and the trace annihilates $D\Qf$. The Leibniz rule makes the
pairing compatible with $D$. In a finite-dimensional complex with such
a perfect pairing, the boundaries are the orthogonal complement of the
cycles. Inclusion follows from the Leibniz rule. For equality, compatibility
with the perfect pairing identifies the maps
$D:\Qf^{52-k}\to \Qf^{53-k}$ and $D:\Qf^k\to \Qf^{k+1}$ as adjoints up to
sign; their ranks agree. Rank--nullity then gives equality with the
orthogonal complement of the cycles in each degree. The induced cohomology pairing is therefore perfect. The top monomial
is closed and cannot be a boundary, which proves the asserted formal
dimension.
\end{proof}

The finite number of generators of $\N$ and its finite-dimensional
cohomology show that the rational space $N$ is elliptic.

\subsection{The orientation and rational category}
The syzygies in \eqref{eq:syzygies} give the degree-$53$ cocycle
\begin{equation}\label{eq:N-orientation}
 \omega=axz^2\sigma_0\sigma_1
        =ax^5z^2b_1b_2-ax^3z^3b_0b_2+axz^4b_0b_1.
\end{equation}
It maps to the top monomial of $\Qf$, so it represents
a nonzero orientation class.

For subsequent calculations introduce the monomials
\begin{equation}\label{eq:T-and-L}
 \begin{gathered}
 T_1=axz^4b_0b_1,\quad T_2=ax^3z^3b_0b_2,\quad
 T_3=ax^5z^2b_1b_2,\\
 T_4=ax^5z^2b_0b_3,\quad T_5=x^3z^6b_0,\quad
 L_1=x^5z^5b_1,\quad L_2=x^7z^4b_2.
 \end{gathered}
\end{equation}
Thus $\omega=T_3-T_2+T_1$. Define
\begin{equation}\label{eq:N-homotopy}
 \eta=-\tfrac12xz^4b_0b_2+3x^3z^3b_1b_2,
 \qquad \zeta=\omega+D\eta.
\end{equation}
Two applications of the Leibniz rule give
\begin{align*}
 D(xz^4b_0b_2)&=4T_2+L_2-T_5+2T_1,\\
 D(x^3z^3b_1b_2)&=3T_3+L_2+T_2-L_1.
\end{align*}
Substitution in \eqref{eq:N-homotopy} yields
\begin{equation}\label{eq:N-long}
 \zeta=10T_3+\tfrac52L_2+\tfrac12T_5-3L_1\in\m_{\N}^{10}.
\end{equation}
The four monomials in \eqref{eq:N-long} have lengths $10,12,10,11$.

For an ordered monomial $m$, write $[m]^\vee$ for its coefficient
functional. The pullback of the top-coefficient trace of $\Qf$ is
\[
 \lambda_{\mathrm{can}}=[T_3]^\vee+3[L_1]^\vee.
\]
Indeed, the only additional word reducing to the top monomial is
$x^5z^5b_1\equiv3T_3\pmod{I_{\mathrm{reg}}}$. The coefficient identity
\begin{equation*}
 D^*[ax^5z^5b_0]^\vee
   =[T_1]^\vee+[T_2]^\vee+[T_4]^\vee
      +6[T_5]^\vee+[L_1]^\vee
\end{equation*}
follows from the differential on the generators. Therefore the functional
\begin{equation}\label{eq:N-trace}
 \lambda=\lambda_{\mathrm{can}}-3D^*[ax^5z^5b_0]^\vee
   =[T_3]^\vee-3[T_1]^\vee-3[T_2]^\vee
      -3[T_4]^\vee-18[T_5]^\vee
\end{equation}
is closed, satisfies $\lambda(\omega)=\lambda(\zeta)=1$, and has support
of length at most ten. Lemma~\ref{lem:toomer-trace} proves
\begin{equation}\label{eq:catN}
 e_0(N)=\cat_0(N)=10.
\end{equation}

\subsection{The obstruction to nilpotency length ten}
\begin{proposition}\label{prop:hnN}
The model $\N$ satisfies $\hn(\N)=11$.
\end{proposition}
\begin{proof}
The upper bound follows from Proposition~\ref{prop:finite-quotient}.
For the lower bound we use the cocycle $\omega$ and the monomials
of \eqref{eq:T-and-L}.

Suppose an acyclic differential ideal $I$ contains $\m_{\N}^{11}$.
There is a functional $\ell:\N^{53}\to\Q$ such that
\begin{equation}\label{eq:ell}
 \ell(I^{53})=0,\qquad \ell D=0,\qquad\ell(\omega)=1.
\end{equation}
Indeed, $\omega$ cannot lie in $I^{53}+D(\N^{52})$: an equality
$\omega=i+Dy$ would make $i$ a cycle in $I$, and its primitive in $I$
would make $\omega$ exact. Choose a linear functional on the quotient
by that subspace and normalize it on $\omega$.
Since $I$ is an ideal, \eqref{eq:ell} also gives
\begin{equation}\label{eq:ell-pair}
 \ell(i s)=0\qquad(i\in I,\ |i|+|s|=53).
\end{equation}

Put
\[
 L_3=x^9z^3b_3,\qquad L_4=ax^7zb_1b_3,\qquad
 L_5=ax^9b_2b_3.
\]
The lengths of $L_1,\ldots,L_5$, with $L_1,L_2$ as in
\eqref{eq:T-and-L}, are $11,12,13,11,12$. Every $L_i$ therefore lies in
$\m_{\N}^{11}$.
The four identities
\begin{equation*}
 \begin{aligned}
 D(axzb_0b_1b_3)&=T_4-T_1-L_4,\\
 D(ax^3b_0b_2b_3)&=T_4-T_2-L_5,\\
 D(x^5z^2b_1b_3)&=T_4+3T_3+2L_4+L_3-L_1,\\
 D(xz^4b_0b_2)&=4T_2+2T_1-T_5+L_2
 \end{aligned}
\end{equation*}
follow by the Leibniz rule. Applying $\ell$ to the first three identities
gives $\ell(T_1)=\ell(T_2)=\ell(T_4)$ and
$\ell(T_4)=-3\ell(T_3)$. The normalization of $\omega$ fixes
$\ell(T_3)=1$, and the fourth identity then gives $\ell(T_5)=-18$.
Thus
\begin{equation}\label{eq:forced-trace}
 \ell(T_3)=1,\qquad
 \ell(T_1)=\ell(T_2)=\ell(T_4)=-3,\qquad
 \ell(T_5)=-18.
\end{equation}
These forced values agree with the coefficients of the closed functional
$\lambda$ in \eqref{eq:N-trace}.

Since $x^{11}\in I$ is closed, acyclicity supplies $\pi\in I^{21}$ with
$D\pi=x^{11}$. We determine all possible primitives in the ambient algebra.
The degree-$21$ monomials are
\[
 ax^9,\ ax^6z,\ ax^3z^2,\ az^3,
\]
and
\[
 x^5b_0,\ x^2zb_0,\ x^4b_1,\ xzb_1,\ x^3b_2,\ zb_2,\ x^2b_3.
\]
Comparing the parts without $a$ in the equation $D\pi=x^{11}$ gives
\[
 \pi=x^5b_0+\mu_0\nu_0+\mu_1\nu_1+\mu_2\nu_2+aP_a,
\]
where
\[
 \begin{gathered}
 \nu_0=x^4b_1-x^2zb_0,\qquad \nu_1=x^3b_2-xzb_1,\qquad
 \nu_2=x^2b_3-zb_2,\\
 \mu_i\in\Q,\qquad P_a\in\Span_{\Q}\{x^9,x^6z,x^3z^2,z^3\}.
 \end{gathered}
\]
One has $D\nu_0=0$, $D\nu_1=ax\sigma_0$, and $D\nu_2=2a\sigma_1$.
The last two images are linearly independent, and $D(aP_a)=0$.
Therefore $\mu_1=\mu_2=0$, and all primitives are exactly
\begin{equation}\label{eq:all-primitives}
 \pi=x^5b_0+\mu_0(x^4b_1-x^2zb_0)+aP_a.
\end{equation}

Multiply \eqref{eq:all-primitives} by $xz^5$, of degree $32$.
The terms $x^6z^5b_0$ and $\mu_0x^5z^5b_1$, and the products from the
first three monomials of $aP_a$, lie in $\m_{\N}^{11}$. The remaining product from
$aP_a$ is a multiple of
\[
 axz^8=-D(axz^5b_3).
\]
The only term detected by $\ell$ is therefore $-\mu_0T_5$.
Equations~\eqref{eq:ell-pair} and \eqref{eq:forced-trace} imply
\[
 0=\ell(\pi xz^5)=-\mu_0\ell(T_5)=18\mu_0,
\]
so $\mu_0=0$.

Now multiply the primitive by $az^2b_3$, again of degree $32$.
The term $aP_a$ disappears because $a^2=0$, and moving $a$ past $b_0$
produces a minus sign. Hence
\[
 0=\ell(\pi az^2b_3)=-\ell(ax^5z^2b_0b_3)
                =-\ell(T_4)=3,
\]
a contradiction. No acyclic differential ideal contains
$\m_{\N}^{11}$. Together with the regular quotient of nilpotency length
eleven, this proves $\hn(\N)=11$.
\end{proof}

\section{Smooth realization}\label{sec:smooth}
The finite quotient supplies the Poincaré duality required for rational
surgery. In the present dimension the realization criterion has no
middle-dimensional arithmetic obstruction.

\begin{theorem}[Barge--Sullivan, the case $4\nmid d$]\label{thm:surgery}
Let $X$ be a simply connected rational space of finite type whose
cohomology is a Poincaré duality algebra of formal dimension $d\geq5$.
If $4\nmid d$, there is an oriented, simply connected, closed smooth
$d$-manifold $V$ such that $V_{\Q}\simeq X$. The realization can be
chosen with all rational Pontryagin classes zero.
\end{theorem}
\begin{proof}
Prescribe the rational model of $X$ and set every rational Pontryagin
class to zero. For $d\geq5$ with $4\nmid d$, Sullivan's realization theorem
provides a simply connected, closed smooth manifold with these prescribed
data \mbox{\cite[\mbox{Theorem}~13.2(i), p.~321]{Sullivan}}.
In particular, its minimal rational model is that of $X$, so
$V_{\Q}\simeq X$. Both the simple connectivity of $V$ and the preservation
of the full rational homotopy type follow from Sullivan's theorem.
No further arithmetic condition remains in these dimensions.

Barge gives the corresponding rational surgery formulation
\mbox{\cite[\mbox{Theorem}~1(a), p.~471]{Barge}}: a closed oriented manifold and
a rational homology equivalence to $X$, with prescribed rational stable
normal Pontryagin classes. His definition does not require the realizing
manifold to be simply connected \mbox{\cite[p.~470]{Barge}}; that additional
conclusion is supplied here by Sullivan, not by Barge's formulation alone.
Vanishing of the stable normal Pontryagin classes is equivalent to
vanishing of the tangent Pontryagin classes, since the two total classes
are multiplicative inverses.
\end{proof}

\Needspace*{12\baselineskip}
\begin{proposition}\label{prop:smoothN}
There is an oriented, simply connected, closed smooth manifold $M^{53}$
whose minimal rational model is $\N$. It is rationally elliptic and can
be chosen with $p_i(TM)=0$ in rational cohomology for all $i\geq1$.
Moreover,
\[
 e_0(M)=\cat_0(M)=\Cat_0(M)=\nilh(M)=10,
 \qquad \hn(M)=11.
\]
\end{proposition}
\begin{proof}
The minimal algebra $\N$ is simply connected, has seven generators,
and is quasi-isomorphic to the finite-dimensional quotient $\Qf$.
Proposition~\ref{prop:finite-quotient} proves Poincaré duality of formal
dimension $53$, with orientation represented by $\omega$ in
\eqref{eq:N-orientation}. Its Sullivan realization $N$ is a simply
connected rational space of finite type. Since $53\geq5$ and
$53\equiv1\pmod4$, Theorem~\ref{thm:surgery} provides a manifold
$M^{53}$ with $M_{\Q}\simeq N$ and the stated Pontryagin classes.

An isomorphism of minimal models identifies augmentation-ideal powers
and graded acyclic differential ideals. Hence homology nilpotency depends
only on the rational homotopy type, and Proposition~\ref{prop:hnN} gives
$\hn(M)=11$. Equation~\eqref{eq:catN} gives $\cat_0(M)=e_0(M)=10$.
The Poincaré duality comparisons in \eqref{eq:comparisons}, using
\mbox{\cite[\mbox{Theorem}~1.1]{CFL}}, identify strong category and
homotopical nil-length with these values. The dual of the generating space of $\N$ identifies with
$\pi_*(M)\otimes\Q$ \mbox{\cite{FHT}}, so
$\dim_{\Q}\bigl(\pi_*(M)\otimes\Q\bigr)=7$. Finite-dimensional
rational homotopy and cohomology are preserved by rational equivalence,
so $M$ is rationally elliptic.
\end{proof}

The product $M\times S^2$ is a simply connected, closed smooth manifold
of dimension $55$. Its rational model is $\A=\N\otimes\Lambda(u,v)$.
Thus determining $\hn(\A)$ determines the invariant of this ordinary
product of manifolds; no separate realization of the product is needed.
The construction specifies a rational homotopy type, rather than a
diffeomorphism type or an integral category value.

\section{The product and its orientation functionals}\label{sec:product}
We now fix
\begin{equation}\label{eq:product-model}
 \A=\N\otimes\Lambda(u,v),\qquad
 |u|=2,\quad |v|=3,\quad Du=0,\quad Dv=u^2.
\end{equation}
The order of the odd letters is $a,b_0,b_1,b_2,b_3,v$. We first establish
the lower bound and then choose the orientation functional used for the
upper bound.

\subsection{The orientation gives the lower bound}
The regular eliminations for $\N$, followed by $(v,u^2)$, give a
quasi-isomorphism
\begin{equation*}
 \A\longrightarrow\overline{\Qf}
     =\Qf\otimes\Lambda(u)/(u^2).
\end{equation*}
The quotient has a perfect top-coefficient pairing in degree $55$.
Its top monomial is $ax^5z^2b_1b_2u$, and
$\overline{\Qf}^{54}=0$. In particular, the cocycle
$\Omega=\omega u$ maps to a nonzero orientation and is not exact.

Using \eqref{eq:N-homotopy} and \eqref{eq:N-long}, set
\begin{equation}\label{eq:Xi}
 \Xi=\zeta u=\Omega+D(\eta u)
   =\left(10T_3+\tfrac52L_2+\tfrac12T_5-3L_1\right)u.
\end{equation}
The four terms in \eqref{eq:Xi} have lengths $11,13,11,12$. Hence $\Xi$ is a nonexact
cocycle in $\m_{\A}^{11}$. An acyclic ideal containing
$\m_{\A}^{11}$ would make $\Xi$ exact, which proves
\begin{equation}\label{eq:product-lower}
 \hn(\A)\geq11.
\end{equation}
This uses the long representative $\Xi$; the three terms of $\Omega$
have lengths $11,10,9$. The category product theorem
\eqref{eq:cat-product} also gives
\begin{equation}\label{eq:cat-product-value}
 \cat_0(N\times S^2)=\cat_0(M\times S^2)=11.
\end{equation}

\subsection{An optimal ideal must mix sphere weights}
Give $\A$ the sphere-weight grading
\begin{equation}\label{eq:sphere-weight}
 \mathrm{wt}_S(\N)=0,\qquad \mathrm{wt}_S(u)=1,\qquad
 \mathrm{wt}_S(v)=2.
\end{equation}
The differential preserves this grading. It is independent of ordinary
cohomological degree and word length.

\Needspace*{9\baselineskip}
\begin{proposition}\label{prop:weight-obstruction}
There is no acyclic differential ideal of $\A$ that contains
$\m_{\A}^{12}$ and is homogeneous in sphere weight.
\end{proposition}
\begin{proof}
Suppose $I$ were such an ideal. The sphere-weight-one part of $\A$ is
exactly $\N u$, so $I_1=I_Nu$ for a graded subspace $I_N\subset\N$.
Multiplication by $\N$ and the differential make $I_N$ a differential
ideal. Since $I$ is a direct sum of its weight subcomplexes, its
acyclicity implies that of $I_1$, hence of $I_N$. But
$\m_{\N}^{11}u\subset\m_{\A}^{12}\subset I$, so
$\m_{\N}^{11}\subset I_N$, contrary to Proposition~\ref{prop:hnN}.
\end{proof}

The saving of one unit, from the tensor bound twelve to eleven, therefore
requires an ideal that mixes sphere weights. This is an additional
constraint on the quotient construction, not a change in the
cohomological grading used to define $\hn$.

\subsection{Two exact corrections of the orientation functional}
The functional
\[
 \eps_1=\lambda\otimes[u]^\vee
\]
is closed, has support of length at most eleven, and satisfies
$\eps_1(\Omega)=1$. Nevertheless, it cannot give the required acyclic
ideal through a weight-homogeneous contraction construction, by
Proposition~\ref{prop:weight-obstruction}. There is also a direct
obstruction in its radical. The sphere-weight-one part of any primitive
$y$ of $x^{11}u$ is $u\pi$, with $D\pi=x^{11}$. Formula
\eqref{eq:all-primitives} and the coefficients of \eqref{eq:N-trace}
give
\begin{equation}\label{eq:uncorrected-obstruction}
 \eps_1(yaz^2b_3)=\lambda(\pi az^2b_3)=3.
\end{equation}
No primitive of this length-twelve cycle therefore belongs to $O_{\eps_1}$.
The obstruction is to multiplicatively compatible primitives, rather
than to the support length of the orientation functional alone.

Define two functionals of support degree $56$ by
\[
 \tau_0=[ax^7z^4b_2]^\vee-3[x^7z^7]^\vee,
 \qquad \tau_3=[ax^2z^5b_1u^3]^\vee.
\]
For $\alpha,\beta\in\Q$, put
\begin{equation}\label{eq:trace-family}
 \eps_{\alpha,\beta}=\eps_1+\alpha D^*\tau_0+\beta D^*\tau_3.
\end{equation}
The two exact terms have sphere weights zero and three. Both are needed;
Section~\ref{sec:mechanism} will characterize their roles separately.

To describe \eqref{eq:trace-family} in coefficients, let
\begin{equation}\label{eq:pure-data}
 \begin{gathered}
 R=\Q[x,z,u],\qquad c_i=b_i\ (0\leq i\leq3),\qquad c_4=v,\\
 \K=R\otimes\Lambda(c_0,c_1,c_2,c_3,c_4),\qquad
 (f_0,f_1,f_2,f_3,f_4)=(x^6,x^4z,x^2z^2,z^3,u^2),\\
 \partial R=0,\qquad\partial c_i=f_i,\qquad
 \q=(f_0,f_1,f_2,f_3,f_4).
 \end{gathered}
\end{equation}
Extend $\theta$ by $\theta u=\theta v=0$. Then
$\A=\Lambda(a)\otimes\K$ and $D=\partial+a\theta$ on $\K$.
The subscript in $\K_j$ denotes the number of letters $c_i$, called
Koszul degree. The differential $\partial$ lowers this degree by one,
whereas $\partial^*$ raises it by one.

For a functional $\xi$ on $\K$, define
\[
 (\ah\xi)(as)=\xi(s),\qquad (\ah\xi)(s)=0\quad(s\in\K).
\]
The operation $\ah$ adds $a$ to the supports; it is neither multiplication
by $a$ nor contraction. For $\xi\in(\K_j)^\vee$, right contraction gives
$a\ctr(\ah\xi)=(-1)^j\xi$. Introduce the coefficient notation
\[
 e_{pqr}=[x^pz^qu^r]^\vee,\quad
 e^i_{pqr}=[x^pz^qu^rc_i]^\vee,\quad
 e^{ij}_{pqr}=[x^pz^qu^rc_ic_j]^\vee\quad(i<j).
\]
The superscripts $i,j$ name the odd letters, not products of functionals.

\Needspace*{15\baselineskip}
We can write
\begin{equation}\label{eq:trace-split}
 \eps_{\alpha,\beta}=\Phi_{\alpha,\beta}+\ah\Psi_{\alpha,\beta},
\end{equation}
where
\begin{equation}\label{eq:Phi}
 \begin{aligned}
 \Phi_{\alpha,\beta}&=-3\alpha e^0_{170}-18e^0_{361}
       +6\beta e^1_{063}-3\alpha e^1_{360}\\
       &+2\beta e^2_{253}+2\alpha e^2_{550},
 \end{aligned}
\end{equation}
and
\begin{equation}\label{eq:Psi}
 \begin{aligned}
 \Psi_{\alpha,\beta}&=-3e^{01}_{141}-\alpha e^{02}_{140}
       -3e^{02}_{331}-3e^{03}_{521}\\
   &+\beta e^{12}_{033}-\alpha e^{12}_{330}+e^{12}_{521}
       +\beta e^{13}_{223}\\
   &+\beta e^{14}_{251}+\alpha e^{23}_{710}.
 \end{aligned}
\end{equation}
In particular, $\Phi$ is supported on $\K_1$ in cohomological degree
$55$, and $\Psi$ on $\K_2$ in degree $52$.

\begin{lemma}\label{lem:trace-properties}
For every $\alpha,\beta\in\Q$, the functional
$\eps_{\alpha,\beta}$ is closed, has support of length at most eleven,
and takes value one on $\Omega$.
\end{lemma}
\begin{proof}
Closedness and the orientation value follow from
\eqref{eq:trace-family}, since $D\Omega=0$ and
$\eps_1(\Omega)=\lambda(\omega)=1$. Expanding the two exact terms by
the Leibniz rule gives \eqref{eq:Phi}--\eqref{eq:Psi}.
For $D^*\tau_0$, the contributions to $[x^7z^4b_3]^\vee$ cancel and
those to $[x^5z^5b_2]^\vee$ have combined coefficient $5-3=2$.
The remaining coefficients are the terms with $\alpha$ in the two
formulas. The five contributions of $D^*\tau_3$ are exactly the terms
with $\beta$. The other five coefficients come from
$\lambda\otimes[u]^\vee$. Each displayed monomial in $\Phi$ has length
at most eleven, and each monomial in $\ah\Psi$ also has length at most
eleven. This proves the support assertion.
\end{proof}

We will construct a quotient for all $\alpha\beta\ne0$. The form of
\eqref{eq:trace-split} permits the same acyclicity argument for every
choice of these two nonzero parameters.

\section{Recognition by polynomial residues}\label{sec:recognition}
The upper bound will follow from a criterion for extensions of pure
Koszul algebras by one closed odd generator. We state it before verifying
its hypotheses for the particular monomials of \eqref{eq:pure-data}.

Let $R$ be a polynomial algebra over $\Q$ on finitely many generators of
positive even degrees, and let
\begin{equation}\label{eq:general-pure}
 \K=\bigl(R\otimes\Lambda(c_1,\ldots,c_s),\partial\bigr),\qquad
 \partial R=0,\qquad\partial c_i=f_i\in R,
\end{equation}
where every $c_i$ is odd. The subscript in $\K_j$ denotes Koszul degree.
Let $a$ be odd of degree at least three, and let $\theta$ be a derivation
of degree $1-|a|$ that preserves Koszul degree and commutes with
$\partial$. Set
\begin{equation}\label{eq:general-extension}
 \A=\Lambda(a)\otimes\K,\qquad
 D=\partial+a\theta\text{ on }\K,\qquad Da=0.
\end{equation}
We assume that this is a simply connected minimal Sullivan algebra.
Write $\q=(f_1,\ldots,f_s)$ and
\[
 B=\partial^*:R^\vee\longrightarrow(\K_1)^\vee,
 \qquad G=(\theta|_R)^*:R^\vee\longrightarrow R^\vee.
\]
The operator $B$ is $R$-linear for contraction; its $i$th component is
$f_i\ctr\xi$ on an input $\xi$. The commutation
$\partial\theta=\theta\partial$ implies $\theta\q\subset\q$.

Suppose $\eps=\Phi+\ah\Psi$ is a homogeneous closed functional on
$\A$, with $\Phi\in(\K_1)^\vee$ and $\Psi\in(\K_2)^\vee$.
The operation $\ah$ and all contractions use the conventions of
Section~\ref{sec:product}.

\begin{definition}\label{def:residue}
The \emph{polynomial residue ideal} of $\eps$, relative to
\eqref{eq:general-extension}, is
\begin{equation}\label{eq:residue}
 \res_\eps=\{\partial\kappa:\kappa\in\K_1,\ 
   \text{there exists }\eta\in\K_2
   \text{ with }\kappa+a\eta\in O_\eps\}\subset R.
\end{equation}
\end{definition}
This is a graded ideal, not merely a set of polynomials. Indeed,
$O_\eps\cap(\K_1+a\K_2)$ is an $R$-module, projection onto $\K_1$ is
$R$-linear, and $\partial:\K_1\to R$ is $R$-linear. All these maps
respect the ordinary grading, with the differential shift where
appropriate. For an ideal $\mathfrak b\subset R$, we use
$\mathfrak b^\perp$ for its annihilator in $R^\vee$.

The element $\kappa+a\eta$ in \eqref{eq:residue} need not be a primitive
of $\partial\kappa$ for the full differential $D$. Its differential may
have an additional $a$-component. The residue records only the polynomial
part, and this is the part detected by scalar dual boundaries. This
distinction permits the criterion to use compatible chains of corrections,
not only primitives of closed polynomials.

For a finite-dimensional graded $R$-submodule $\W\subset R^\vee$, set
\begin{equation}\label{eq:general-E}
 \begin{gathered}
 \F(\W)=\W\oplus\ah\W\oplus B(\W)\oplus\ah B(\W),\\
 \E=\F(\W)+\A\ctr\eps,\qquad J=\E^\perp.
 \end{gathered}
\end{equation}

\Needspace*{15\baselineskip}
\begin{theorem}[Residue recognition]\label{thm:recognition}
Suppose \eqref{eq:general-pure}--\eqref{eq:general-extension} and
$\eps=\Phi+\ah\Psi$ satisfy the following conditions.
\begin{enumerate}
\item $H(\K)$ is a finite-dimensional Poincaré duality algebra whose
orientation has Koszul degree two, and $\Psi$ is a closed pure
orientation functional. In particular, $H_j(\K)=0$ for $j>2$.
\item A finite-dimensional graded $R$-submodule $\W\subset R^\vee$
contains $(R/\q)^\vee$ and every scalar coefficient functional of
$\Phi$ and $\Psi$. It is stable under $G$.
\item $\res_\eps^\perp\subset\W$.
\end{enumerate}
Then $\E\hookrightarrow\A^\vee$ is a quasi-isomorphism, and $J$ is an
acyclic differential ideal.
If $\eps$ and all functionals $\ah\xi$, $\xi\in\W$, have support of
length at most $n$, then $\m_{\A}^{n+1}\subset J$ and $\hn(\A)\leq n$.
\end{theorem}

Here $\A\ctr\eps$ means the span of all contractions of $\eps$.
The scalar coefficients in condition~(2) are defined by expansion in
the exterior basis of the $c_i$: the coefficients of
$\Phi\in(\K_1)^\vee$ and $\Psi\in(\K_2)^\vee$ are elements of
$R^\vee$. The theorem retains both ordinary cohomological degree and
Koszul degree. The vanishing in condition~(1) follows from the perfect
pairing between Koszul degrees $j$ and $2-j$, since negative Koszul
degrees are absent. The theorem requires $\Psi$ to orient the pure
cohomology; it does not impose a separate orientation normalization on
$\eps$ for $\A$. No homogeneity of $J$ in any additional weight is imposed.

\begin{corollary}\label{cor:residue-tail}
Let $R=\Q[y_1,\ldots,y_m]$, and let $\W_{n-1}\subset R^\vee$
consist of the coefficient functionals supported on polynomial monomials
of length at most $n-1$. If conditions~(1) and~(2) of
Theorem~\ref{thm:recognition} hold with $\W=\W_{n-1}$, and
\[
 (y_1,\ldots,y_m)^n\subset\res_\eps,
\]
then condition~(3) holds. If $\eps$ has support of length at most $n$,
the ideal in \eqref{eq:general-E} is acyclic and contains
$\m_{\A}^{n+1}$.
\end{corollary}
\begin{proof}
A functional annihilating $\res_\eps$ annihilates every polynomial
monomial of length at least $n$, so it belongs to $\W_{n-1}$.
The generators $\ah\W_{n-1}$ have support of length at most $n$.
The conclusion is Theorem~\ref{thm:recognition}.
\end{proof}

The first condition is a statement about the pure algebra $\K$, whose
Poincaré pairing will be obtained from a regular quotient. The residue
condition is an inclusion in the ordinary polynomial ring. The
construction will be applied after both conditions have been verified.

\section{Proof of the recognition theorem}\label{sec:recognition-proof}
The proof proceeds through a pure dual subcomplex. Residues first give
primitives for all ambient boundaries that meet this subcomplex. Pure
duality then controls the additional terms involving $\Phi$, and
contraction by $a$ identifies the extension that must be compared with
$\A^\vee$.

\subsection{The polynomial module and the image of contraction}
For a scalar functional $\xi\in R^\vee$, direct transposition of
\eqref{eq:general-extension} gives
\begin{equation}\label{eq:poly-D}
 \begin{gathered}
 D^*\xi=B\xi,\qquad D^*(\ah\xi)=G\xi-\ah B\xi,\\
 D^*B\xi=0,\qquad D^*(\ah B\xi)=BG\xi.
 \end{gathered}
\end{equation}
The minus sign comes from $D(as)=-a\partial s$ for $s\in\K$.
For the last identity use $\theta^*B=BG$, a consequence of
$\partial\theta=\theta\partial$.

The four summands defining $\F(\W)$ in \eqref{eq:general-E} have
separate odd-letter types and form a direct sum. They are stable under
$D^*$ by \eqref{eq:poly-D} and the $G$-stability of $\W$.
Polynomial contractions preserve them because $\W$ is an $R$-module
and $B$ is $R$-linear. Contraction by $c_i$ sends $B\xi$ or
$\ah B\xi$ to $f_i\ctr\xi$ or its $\ah$ copy, up to sign;
contraction by $a$ removes $\ah$. This proves stability of $\F(\W)$ under
all contractions. It is the smallest differential contraction module
containing $\ah\W$: contraction by $a$ gives $\W$, the differential
gives $B(\W)$, and \eqref{eq:poly-D} then gives $\ah B(\W)$.
Since $D^*\eps=0$, formula~\eqref{eq:contraction-differential} makes
$\A\ctr\eps$ differential as well. We obtain a finite-dimensional
differential contraction submodule $\E\subset\A^\vee$.

Put
\begin{equation}\label{eq:pure-subcomplex}
 \begin{gathered}
 I_\Psi=\Ann_R(\Psi),\qquad
 V_1=\sum_{i=1}^{s}R\ctr(c_i\ctr\Psi),\qquad V_2=R\ctr\Psi,\\
 \Pmod=\W\oplus\bigl(B(\W)+V_1\bigr)\oplus V_2\subset\K^\vee.
 \end{gathered}
\end{equation}
The direct summands in \eqref{eq:pure-subcomplex} have Koszul degrees
zero, one, and two.
All polynomial contractions of the coefficients of $\Psi$ belong to
$\W$, by condition~(2). In particular, contracting $\Psi$ by two
letters $c_i$ adds no scalar functionals outside $\W$.

\begin{lemma}\label{lem:a-image}
The image of $a\ctr:\E\to\E$ is $\Pmod$. The restriction of $D^*$
to $\Pmod$ is $\partial^*$, so $\Pmod$ is a pure subcomplex of
$\K^\vee$.
\end{lemma}
\begin{proof}
Contraction of $\F(\W)$ by $a$ has image $\W+B(\W)$.
Also $a\ctr\eps=\Psi$, since $\Psi$ has Koszul degree two.
Contraction of $\A\ctr\eps$ by $a$ therefore gives all contractions
of $\Psi$: these have components $V_2$, $V_1$, and scalar components
in $\W$. This proves the image assertion, including the inclusion
$\Pmod\subset\E$.

Right contraction by $a$ commutes with $D^*$: for every functional
$\xi$ and every $s$, both composites evaluate as $\xi((Ds)a)$ because
$Da=0$. Its image is a subcomplex. A functional supported on $\K$ and
zero on $a\K$ sees only the pure part of $D$, so its differential is
$\partial^*$.
\end{proof}

\subsection{Residues control ambient boundaries}
The next lemma treats both the pure contractions $V_1$ and the terms
$I_\Psi\ctr\Phi$. The stronger statement is what will make the same
argument control the extension by $a$.

\begin{lemma}[Saturation of boundaries]\label{lem:saturation}
Under condition~(3) of Theorem~\ref{thm:recognition},
\begin{equation}\label{eq:saturation}
 \bigl(V_1+I_\Psi\ctr\Phi\bigr)\cap B(R^\vee)\subset B(\W).
\end{equation}
In particular, every element of $V_1$ that is a boundary in the full
pure dual has a primitive in $\W$.
\end{lemma}
\begin{proof}
Let $B(h)\in V_1+I_\Psi\ctr\Phi$, with $h\in R^\vee$, and take
$\rho=\kappa+a\eta\in O_\eps$ with $\kappa\in\K_1$ and
$\eta\in\K_2$. The two odd-letter types of $\rho\ctr\eps=0$ give
\begin{equation}\label{eq:orphan-components}
 \kappa\ctr\Psi=0,\qquad
 \kappa\ctr\Phi=-\eta\ctr\Psi.
\end{equation}
The first identity is the component with $a$ and one letter $c_i$;
the second is the scalar component. Every element of $V_1$ vanishes on
$\kappa$, because its value is, up to the right-contraction sign, an
evaluation of $\kappa\ctr\Psi$ on a polynomial multiple of a letter
$c_i$.

For a polynomial $P\in I_\Psi$, the second identity in
\eqref{eq:orphan-components} gives
\[
 (P\ctr\Phi)(\kappa)=(\kappa\ctr\Phi)(P)
                    =-(\eta\ctr\Psi)(P)=0,
\]
since $P\ctr\Psi=0$. Therefore $B(h)(\kappa)=0$, or equivalently
\[
 h(\partial\kappa)=0.
\]
This holds for every admissible $\kappa+a\eta$ in the definition of
$\res_\eps$. It follows that $h\in\res_\eps^\perp\subset\W$, proving
\eqref{eq:saturation}. The last assertion follows by setting the
$I_\Psi\ctr\Phi$ component to zero.
\end{proof}

The residue ideal is used as an ideal throughout this proof. Multiples
of the admissible elements are themselves admissible, so the conclusion
controls arbitrary polynomial combinations, rather than only the
particular residues used to verify condition~(3).

\subsection{Comparison with the pure dual}
Write $H^{[j]}$ for cohomology in the increasing Koszul degree on a pure
dual complex. This bracketed index is distinct from the ordinary
support degree. The duality map \eqref{eq:PD-contraction}, applied to $\Psi$, gives
\begin{equation}\label{eq:pure-duality}
 H_j(\K)\xrightarrow{\ \cong\ }H^{[2-j]}(\K^\vee),
 \qquad [\kappa]\longmapsto[\kappa\ctr\Psi],
 \qquad 0\leq j\leq2,
\end{equation}
with the ordinary degree shifts understood. The orientation degree is
even, so on an element of $\K_j$ the compatibility formula is
$\partial^*(\kappa\ctr\Psi)=(-1)^j(\partial\kappa)\ctr\Psi$.
These signs do not affect the induced duality isomorphisms.

\Needspace*{7\baselineskip}
\begin{lemma}\label{lem:pure-comparison}
The inclusion $\Pmod\hookrightarrow\K^\vee$ is a quasi-isomorphism.
\end{lemma}
\begin{proof}
First,
\begin{equation}\label{eq:ann-in-q}
 I_\Psi\subset\q.
\end{equation}
Indeed, $P\in I_\Psi$ is a pure cycle and $P\ctr\Psi=0$. The duality
map \eqref{eq:pure-duality} makes its class in
$H_0(\K)=R/\q$ zero, so $P\in\q$.

In Koszul degree zero, the kernel of $B$ in $R^\vee$ is
$(R/\q)^\vee$, entirely contained in $\W$. The inclusion is therefore
an isomorphism on degree-zero cohomology.

In degree one, suppose a cycle $\xi\in B(\W)+V_1$ is a boundary in
$\K^\vee$. Subtracting a component in $B(\W)$ leaves an element of
$V_1\cap B(R^\vee)$. Lemma~\ref{lem:saturation} places this component
in $B(\W)$ as well. The element $\xi$ is already a boundary in $\Pmod$,
proving injectivity. For surjectivity, \eqref{eq:pure-duality} represents
every class of $H^{[1]}(\K^\vee)$ by contraction of $\Psi$ with a cycle
in $\K_1$. Such a contraction belongs to $V_1$.

In degree two, $V_2$ is a cyclic $R$-module isomorphic to $R/I_\Psi$.
The differential of a contraction by $c_i$ is
$\partial^*(c_i\ctr\Psi)=-f_i\ctr\Psi$. Therefore
\[
 \partial^*V_1=\q V_2,
 \qquad H^{[2]}(\Pmod)=V_2/\q V_2\cong R/(I_\Psi+\q)=R/\q,
\]
where the last equality uses \eqref{eq:ann-in-q}. The map to
$H^{[2]}(\K^\vee)$ is the duality isomorphism on $H_0(\K)$ in
\eqref{eq:pure-duality}. Finally, both complexes have zero cohomology
in Koszul degrees above two, by condition~(1) and the definition of
$\Pmod$. This proves the lemma in all degrees.
\end{proof}

\subsection{The annihilator of the pure orientation}
We now use pure duality once more, this time to control those contractions
of $\Phi$ whose corresponding contraction of $\Psi$ vanishes.

\begin{lemma}\label{lem:ann-Phi}
One has
\begin{equation}\label{eq:ann-Phi}
 I_\Psi\ctr\Phi\subset V_1+B(\W).
\end{equation}
\end{lemma}
\begin{proof}
The closedness of $\eps=\Phi+\ah\Psi$ gives
\begin{equation}\label{eq:closed-components}
 \partial^*\Psi=0,\qquad \partial^*\Phi=-\theta^*\Psi.
\end{equation}
For a homogeneous polynomial $P\in I_\Psi$, transposing the Leibniz
rule for the even derivation $\theta$ yields
\[
 P\ctr(\theta^*\Psi)=\theta^*(P\ctr\Psi)+(\theta P)\ctr\Psi
                                      =(\theta P)\ctr\Psi.
\]
Combining this with \eqref{eq:closed-components} gives
$\partial^*(P\ctr\Phi)=-(\theta P)\ctr\Psi$.
By \eqref{eq:ann-in-q} and $\theta\q\subset\q$, there exists
$\kappa_0\in\K_1$ such that $\partial\kappa_0=\theta P$.
The contraction sign in odd Koszul degree gives
\[
 \partial^*(\kappa_0\ctr\Psi)=-(\partial\kappa_0)\ctr\Psi.
\]
It follows that $P\ctr\Phi-\kappa_0\ctr\Psi$ is a pure dual cycle
of Koszul degree one.

The duality \eqref{eq:pure-duality} supplies a cycle
$\kappa_1\in\K_1$ whose contraction represents its cohomology class.
Thus, for some $h\in R^\vee$,
\[
 P\ctr\Phi=(\kappa_0+\kappa_1)\ctr\Psi+B(h).
\]
The first term belongs to $V_1$. The term $B(h)$ therefore belongs to the
left-hand side of \eqref{eq:saturation}, and
Lemma~\ref{lem:saturation} gives $B(h)\in B(\W)$. This proves
\eqref{eq:ann-Phi}.
\end{proof}

The choices of $\kappa_0$, $\kappa_1$, and $h$ in this proof only
establish an inclusion of spaces. No $R$-linear choice of primitives is
required: $V_1+B(\W)$ is already an $R$-module, and the inclusion is
proved for every polynomial in $I_\Psi$.

\subsection{Exact contraction and the comparison of extensions}
\begin{lemma}\label{lem:a-exact}
The kernel and image of $a\ctr:\E\to\E$ are both $\Pmod$.
\end{lemma}
\begin{proof}
Lemma~\ref{lem:a-image} identifies the image. For the kernel, observe
that a contraction of $\eps$ by two letters $c_i$ lies in $\ah\W$,
while a contraction by one such letter belongs to
$\W+\ah V_1$. Contractions involving $a$ lie in $\Pmod$.
It follows that
\begin{equation}\label{eq:E-decomposition}
 \E=\F(\W)+V_1+\ah V_1+V_2+R\ctr\eps.
\end{equation}
Consider an element in the kernel of $a\ctr$. In
\eqref{eq:E-decomposition}, the component of the form $\ah\K_2^\vee$
comes only from $P\ctr\eps$ for a polynomial $P$; its contraction is
$P\ctr\Psi$. It must vanish separately, so $P\in I_\Psi$.
The remaining components with $\ah$ must also cancel. The element is
therefore in $\Pmod+I_\Psi\ctr\Phi$, which is $\Pmod$ by
Lemma~\ref{lem:ann-Phi}. Conversely, every functional in $\Pmod$ is
zero on $a\K$ and is killed by contraction by $a$.
\end{proof}

\begin{proof}[Completion of the proof of Theorem~\ref{thm:recognition}]
In the full dual $\A^\vee$, the kernel and image of contraction by $a$
are $\K^\vee$. Lemma~\ref{lem:a-exact} gives a commutative diagram
\begin{equation}\label{eq:comparison-diagram}
 \begin{array}{ccccccccc}
 0&\longrightarrow&\Pmod&\longrightarrow&\E
      &\xrightarrow{\ a\ctr\ }&\Sigma_a\Pmod&\longrightarrow&0\\
 &&\big\downarrow&&\big\downarrow&&\big\downarrow\\
 0&\longrightarrow&\K^\vee&\longrightarrow&\A^\vee
      &\xrightarrow{\ a\ctr\ }&\Sigma_a\K^\vee&\longrightarrow&0.
 \end{array}
\end{equation}
Here $\Sigma_a$ denotes reindexing of support degrees,
$(\Sigma_a U)^d=U^{d-|a|}$, with unchanged differential. Since right
contraction by $a$ commutes with $D^*$, both rows are short exact
sequences of complexes in this convention. The two outside vertical
maps are quasi-isomorphisms by Lemma~\ref{lem:pure-comparison}.
The comparison of their long exact sequences proves that
$\E\hookrightarrow\A^\vee$ is a quasi-isomorphism.

The differential contraction stability of $\E$ was proved using
\eqref{eq:poly-D}. Lemma~\ref{lem:dual} now shows that $J=\E^\perp$
is an acyclic differential ideal.
For the length assertion, $\E$ is generated as a differential
contraction module by $\eps$ and $\ah\W$. The assumed support bound
is preserved by contraction and by $D^*$, because $D$ is minimal.
Every functional in $\E$ therefore annihilates $\m_{\A}^{n+1}$,
which gives $\m_{\A}^{n+1}\subset J$ and $\hn(\A)\leq n$.
\end{proof}

The comparison \eqref{eq:comparison-diagram} is valid in every ordinary
degree. It requires the pure Poincaré pairing and the saturation
\eqref{eq:saturation}, but no choice of a basis for the action of
$\theta$ on pure cohomology. The residue condition coordinates the
boundary and contraction requirements before the extension by $a$ is
reconstructed.

\section{Applying the criterion to the sphere product}\label{sec:application}
Return to the algebra \eqref{eq:product-model}, its pure part
\eqref{eq:pure-data}, and the family \eqref{eq:trace-family}.
We verify pure duality, the scalar-module conditions, and finally the
residue inclusion. Throughout this section $\alpha\beta\ne0$.

\subsection{Pure Poincaré duality}
\begin{proposition}\label{prop:pure-PD}
The pure algebra $\K$ has Poincaré duality in cohomology of formal
dimension $52$, with orientation of Koszul degree two. Its Koszul
homology vanishes above degree two. For every $\alpha,\beta$, the
functional $\Psi_{\alpha,\beta}$ in \eqref{eq:Psi} is a pure orientation
functional.
\end{proposition}
\begin{proof}
The regular eliminations $(b_0,x^6)$, $(b_3,z^3)$, and $(v,u^2)$ give a
pure quasi-isomorphism
\begin{equation}\label{eq:pure-quotient}
 \K\longrightarrow Q_{\mathrm p}
  =R/(x^6,z^3,u^2)\otimes\Lambda(b_1,b_2),\qquad
 \partial b_1=x^4z,\quad\partial b_2=x^2z^2.
\end{equation}
Lemma~\ref{lem:regular} applies successively because $x^6,z^3,u^2$
are powers of separate polynomial generators and remain regular at the
required stages. The projections preserve Koszul degree. Their kernels
are acyclic in the Koszul complexes as well: the primitive in the proof
of the lemma increases the number of odd generators by one. Since
$Q_{\mathrm p}$ has only two odd generators, $H_j(\K)=0$ for $j>2$.

The quotient \eqref{eq:pure-quotient} has basis
\[
 x^pz^qu^rb_1^eb_2^f,\qquad
 0\leq p\leq5,\quad0\leq q\leq2,\quad0\leq r\leq1,
 \quad e,f\in\{0,1\}.
\]
Its top monomial $x^5z^2b_1b_2u$ has degree $52$.
Extraction of its coefficient gives a perfect cochain pairing by
complementary monomials. Since $Q_{\mathrm p}^{51}=0$, this functional
is closed. The cochain duality argument of Section~\ref{sec:prelim}
therefore gives Poincaré duality in cohomology. Complementary monomials
have a combined total of two odd letters, so the pairing also respects
complementary Koszul degrees.

The pure cycle
\begin{equation}\label{eq:pure-orientation}
 \varpi=xz^2\sigma_0\sigma_1u
       =(x^5z^2b_1b_2-x^3z^3b_0b_2+xz^4b_0b_1)u
\end{equation}
maps to the top monomial of $Q_{\mathrm p}$. Closedness of
$\eps_{\alpha,\beta}$ implies $\partial^*\Psi_{\alpha,\beta}=0$.
The three coefficients of \eqref{eq:Psi} seen by
\eqref{eq:pure-orientation} give
\[
 \Psi_{\alpha,\beta}(\varpi)=1-(-3)+(-3)=1.
\]
The functional $\Psi_{\alpha,\beta}$ therefore induces the nonzero orientation functional
of this pure Poincaré duality algebra.
\end{proof}

\subsection{The scalar module and its length}
Take
\begin{equation}\label{eq:W10}
 \W_{10}=\Span_{\Q}\{e_{pqr}:p+q+r\leq10\}\subset R^\vee.
\end{equation}
It is an $R$-submodule for contraction. A monomial basis for $R/\q$ is
\[
 x^pz^qu^r,\qquad
 0\leq q\leq2,\quad0\leq p\leq5-2q,\quad0\leq r\leq1.
\]
Every such monomial has polynomial length at most six, so
$(R/\q)^\vee\subset\W_{10}$. All scalar coefficient functionals in
\eqref{eq:Phi} and \eqref{eq:Psi} belong to $\W_{10}$ as well.
For all nonnegative triples $(p,q,r)$, the transpose of
$\theta|_R=x^2\partial_z$ satisfies
\begin{equation}\label{eq:G-action}
 G(e_{pqr})=
 \begin{cases}
 (q+1)e_{p-2,q+1,r},&p\geq2,\\
 0,&p\leq1.
 \end{cases}
\end{equation}
The nonzero images in \eqref{eq:G-action} have polynomial support length
one less than the input. Thus $G(\W_{10})\subset\W_{10}$.

Let
\begin{equation}\label{eq:final-ideal}
 \begin{gathered}
 \F=\W_{10}\oplus\ah\W_{10}\oplus B(\W_{10})\oplus\ah B(\W_{10}),\\
 \E_{\alpha,\beta}=\F+\A\ctr\eps_{\alpha,\beta},\qquad
 J_{\alpha,\beta}=\E_{\alpha,\beta}^{\perp}.
 \end{gathered}
\end{equation}
Proposition~\ref{prop:pure-PD} and \eqref{eq:W10} verify the first two
conditions of the recognition theorem. The generators
$\ah\W_{10}$ have support of length at most eleven, as does
$\eps_{\alpha,\beta}$ by Lemma~\ref{lem:trace-properties}.
By Corollary~\ref{cor:residue-tail}, the remaining condition is
\begin{equation}\label{eq:residue-target}
 (x,z,u)^{11}\subset\res_{\eps_{\alpha,\beta}}.
\end{equation}

\subsection{Five admissible elements}
Write $O=O_{\eps_{\alpha,\beta}}$. The sphere relation and the
syzygies \eqref{eq:syzygies} lead to the five elements
\begin{equation}\label{eq:five-elements}
 \begin{aligned}
 \rho_1={}&x^5b_0+\frac3\beta D(v\sigma_0),\\
 \rho_2={}&x^4ub_0+\frac3\beta D(xuvb_1),\\
 \rho_3={}&x^3zb_0-\frac3\beta u^2\sigma_1
                         +\frac6\beta a\sigma_0v,\\
 \rho_4={}&x^3z^3b_1+\frac\alpha\beta D(uvb_1Db_3),\\
 \rho_5={}&\frac u2\bigl(x^3z^3b_2+x^5z^2b_3+6axzb_1b_2\bigr)
                         +\frac3\alpha D(xzb_2b_3).
 \end{aligned}
\end{equation}

\Needspace*{9\baselineskip}
\begin{proposition}\label{prop:five-elements}
The five elements in \eqref{eq:five-elements} lie in
$O\cap(\K_1+a\K_2)$. Their cohomological degrees are
$21,21,23,37,41$. If $\rho_i=\kappa_i+a\eta_i$, their polynomial
residues are
\begin{equation}\label{eq:five-residues}
 \partial\kappa_1=x^{11},\quad
 \partial\kappa_2=x^{10}u,\quad
 \partial\kappa_3=x^9z,\quad
 \partial\kappa_4=x^7z^4,\quad
 \partial\kappa_5=x^5z^5u.
\end{equation}
\end{proposition}
\begin{proof}
The pure part of an exact correction $Dq$ is $\partial(q|_{a=0})$,
so its pure differential is zero. In $\rho_3$ the correction without
$a$ is the pure cycle $-3u^2\sigma_1/\beta$.
In $\rho_5$ the two leading terms each contribute half of $x^5z^5u$.
These observations prove \eqref{eq:five-residues}.

For membership in $O$, both identities in
\eqref{eq:orphan-components} must be checked. The parts without $a$ in
\eqref{eq:five-elements} are
\begin{equation}\label{eq:kappa-parts}
 \begin{aligned}
 \kappa_1={}&x^5b_0+\frac3\beta u^2\sigma_0,\\
 \kappa_2={}&x^4ub_0+\frac3\beta xu^3b_1-\frac3\beta x^5zuv,\\
 \kappa_3={}&x^3zb_0-\frac3\beta u^2\sigma_1,\\
 \kappa_4={}&z^3\left(x^3+\frac\alpha\beta u^3\right)b_1
                         -\frac\alpha\beta x^4z^4uv,\\
 \kappa_5={}&\frac u2\bigl(x^3z^3b_2+x^5z^2b_3\bigr)
                         +\frac3\alpha xz^3\sigma_2.
 \end{aligned}
\end{equation}
Contraction with the ten coefficients in \eqref{eq:Psi} gives
\begin{equation}\label{eq:kappa-Psi}
 \kappa_i\ctr\Psi_{\alpha,\beta}=0\qquad(1\leq i\leq5).
\end{equation}
Here is a componentwise verification of \eqref{eq:kappa-Psi}.
For $\kappa_1$, the only surviving contributions come from the types
$03$ and $13$; their coefficients are $3$ and $-3$ on the same
functional. For $\kappa_2$, the same two types have opposite coefficients;
the term containing $v$ has too many factors $x$ to contribute. In $\kappa_3$, the contributions of types $02$ and $03$
are canceled by those of $12$ and $13$. For $\kappa_4$, only
$-\alpha e^{12}_{330}$ and $\beta e^{12}_{033}$ survive, with opposite
coefficients; the term with $v$ is again excluded by its $x$-exponent.
For $\kappa_5$, the two weight-one contributions of types $02$ and
$03$ cancel the weight-zero contribution of type $02$. Every other
candidate is excluded by one of its polynomial exponents.

The scalar contractions of the parts in \eqref{eq:kappa-parts} are
\begin{equation}\label{eq:kappa-Phi}
 \begin{aligned}
 \kappa_1\ctr\Phi&=0,&\qquad \kappa_2\ctr\Phi&=0,\\
 \kappa_3\ctr\Phi&=-6e_{051},&
 \kappa_4\ctr\Phi&=3\alpha e_{030},\\
 \kappa_5\ctr\Phi&=-6e_{410}-\frac{6\beta}{\alpha}e_{113}.&&
 \end{aligned}
\end{equation}
To verify the compensating contributions of $a\eta_i$, it is enough to
use the five scalar identities
\begin{equation}\label{eq:scalar-contractions}
 \begin{aligned}
 (ax^2b_1v)\ctr\eps&=\beta e_{051},&
 (au^3b_1b_2)\ctr\eps&=\beta e_{030},\\
 (ax^3b_2b_3)\ctr\eps&=\alpha e_{410},&
 (axzb_1b_3)\ctr\eps&=\beta e_{113},\\
 (axzub_1b_2)\ctr\eps&=e_{410}.&&
 \end{aligned}
\end{equation}
Each identity follows from the indicated odd-letter type of
\eqref{eq:Psi}. Terms containing $ab_0v$ or $ab_2v$ have zero
contraction, because these odd-letter types do not occur.

In fact, $a\eta_1=0$, and $a\eta_2$ is a multiple of $ab_0v$.
The part $a\eta_3=(6/\beta)a\sigma_0v$ contributes $6e_{051}$.
The part with $a$ of the exact correction in $\rho_4$ is
\[
 -\frac\alpha\beta a
     \bigl(3u^3b_1b_2+3x^4zub_2v+z^3ub_0v\bigr),
\]
which contributes $-3\alpha e_{030}$. Finally,
\begin{equation}\label{eq:high-homotopy}
 D(xzb_2b_3)=xz^3\sigma_2+ax^3b_2b_3+2axzb_1b_3.
\end{equation}
The $a$-part of $\rho_5$ consequently contributes
$6e_{410}+(6\beta/\alpha)e_{113}$ by
\eqref{eq:scalar-contractions}. These are the negatives of
\eqref{eq:kappa-Phi}. Together with \eqref{eq:kappa-Psi}, they prove
$\rho_i\ctr\eps=0$ for all five elements. The displayed expansions
also show that only the types $\K_1+a\K_2$ occur.
\end{proof}

The differential relations among these elements explain their
coordination. The Leibniz rule and \eqref{eq:syzygies} give
\begin{equation}\label{eq:coherent-corrections}
 \begin{aligned}
 D\rho_1={}&x^{11},\qquad D\rho_2=x^{10}u,\\
 D\rho_3={}&x^9z+a\rho_1,\\
 D\rho_4={}&x^7z^4+a(3x^5z^2b_1+x^3z^3b_0),\\
 D\rho_5={}&u(x^5z^5+ax^7zb_3+4ax^3z^3b_1).
 \end{aligned}
\end{equation}
For the identity involving $D\rho_3$ in \eqref{eq:coherent-corrections},
$D(x^9z)=ax^{11}$ is canceled by
$D(a\rho_1)=-ax^{11}$. The correction is thus compatible with the
preceding primitive. The formula for $D\rho_4$ involves both the sphere relation
and the last differential $Db_3$; the formula for $D\rho_5$ is linked to
the final syzygy by \eqref{eq:high-homotopy}.

The role of the elements \eqref{eq:five-elements} is to test dual
boundaries. Their required membership is in $O_\eps$, not in the ideal
$J_{\alpha,\beta}$ itself. For example,
\[
 B(e_{910})(\rho_3)=e_{910}(\partial\kappa_3)=1,
 \qquad e_{910}\in\W_{10}.
\]
In particular, $\rho_3\notin J_{\alpha,\beta}$. The residue theorem uses all
polynomial multiples of this test element without making it a relation
in the quotient.

\subsection{The residue ideal contains the polynomial tail}
Seven further elements are admissible by their supports alone:
\begin{equation}\label{eq:monomial-tests}
 u^2v,\quad x^3v,\quad z^3b_3,\quad x^6ub_3,\quad
 x^6z^2b_3,\quad x^8b_3,\quad x^6ub_2.
\end{equation}
They belong to $O_\eps\cap\K_1$. For the two terms containing $v$,
the only compatible support in \eqref{eq:Psi} is the type $b_1v$ with
exponent triple $(2,5,1)$; either $u^2$ or $x^3$ prevents divisibility.
For the terms containing $b_3$, the compatible triples are
$(5,2,1)$, $(2,2,3)$, and $(7,1,0)$. Each displayed monomial fails to
divide each of these supports. The term $x^6ub_2$ is excluded from
\eqref{eq:Phi}--\eqref{eq:Psi} by its $x$-exponent or its $u$-exponent.
These divisibility checks prove vanishing under every contraction, as
required for membership in $O_\eps$.

Their residues, together with \eqref{eq:five-residues}, give
\begin{equation}\label{eq:monomial-ideal}
 \begin{aligned}
 \ellid=\bigl(&u^4,x^3u^2,z^6,x^6z^3u,x^6z^5,x^8z^3,x^8z^2u,\\
              &x^{11},x^{10}u,x^9z,x^7z^4,x^5z^5u\bigr)
       \subset\res_{\eps_{\alpha,\beta}}.
 \end{aligned}
\end{equation}
Only the inclusion is needed; no equality of these ideals is assumed.

\Needspace*{8\baselineskip}
\begin{lemma}\label{lem:polynomial-tail}
One has $(x,z,u)^{11}\subset\ellid$.
\end{lemma}
\begin{proof}
Let $x^pz^qu^r$ represent a nonzero monomial in $R/\ellid$.
The generators $u^4$ and $z^6$ force $r\leq3$ and $q\leq5$.
If $r=2$ or $3$, the generator $x^3u^2$ forces $p\leq2$, and
$p+q+r\leq10$.

For $r=0$ or $1$, the remaining generators of
\eqref{eq:monomial-ideal} give the following maximal values of $p$:
\begin{center}
 \begin{tabular}{c@{\qquad}rrrrrr}
 \toprule
 & $q=0$ & $q=1$ & $q=2$ & $q=3$ & $q=4$ & $q=5$\\
 \midrule
 $r=0$ &10&8&8&7&6&5\\
 $r=1$ &9&8&7&5&5&4\\
 \bottomrule
 \end{tabular}
\end{center}
The first row follows respectively from $x^{11},x^9z,x^9z,x^8z^3,
x^7z^4,x^6z^5$. The second follows from $x^{10}u,x^9z,x^8z^2u,
x^6z^3u,x^6z^3u,x^5z^5u$. Each entry satisfies $p+q+r\leq10$.
Every monomial of polynomial length at least eleven is therefore zero
in $R/\ellid$.
\end{proof}

\Needspace*{12\baselineskip}
\begin{theorem}\label{thm:product-acyclicity}
For every $\alpha,\beta\in\Q^\times$, the inclusion
$\E_{\alpha,\beta}\hookrightarrow\A^\vee$ in
\eqref{eq:final-ideal} is a quasi-isomorphism. The ideal
$J_{\alpha,\beta}$ is acyclic and contains $\m_{\A}^{12}$.
Consequently $\hn(N\times S^2)=11$.
\end{theorem}
\begin{proof}
Proposition~\ref{prop:pure-PD} verifies condition~(1) of
Theorem~\ref{thm:recognition}. The scalar module \eqref{eq:W10}
and formula~\eqref{eq:G-action} verify condition~(2).
Proposition~\ref{prop:five-elements}, the tests
\eqref{eq:monomial-tests}, and Lemma~\ref{lem:polynomial-tail} prove
\eqref{eq:residue-target}, so Corollary~\ref{cor:residue-tail} supplies
condition~(3). The support bound is eleven by
Lemma~\ref{lem:trace-properties} and the definition of $\W_{10}$.
The recognition theorem proves the quasi-isomorphism, acyclicity, and
$\m_{\A}^{12}\subset J_{\alpha,\beta}$. This gives
$\hn(\A)\leq11$. The lower bound \eqref{eq:product-lower} gives equality.
\end{proof}

\begin{proof}[Proof of Theorem~\ref{thm:main}]
Proposition~\ref{prop:smoothN} provides a simply connected, closed smooth
manifold $M^{53}$ with the prescribed rational model and zero rational
Pontryagin classes. It also gives $\cat_0(M)=10$ and $\hn(M)=11$.
The product model is \eqref{eq:product-model}, so
Theorem~\ref{thm:product-acyclicity} gives $\hn(M\times S^2)=11$.
Its rational category is given by \eqref{eq:cat-product-value}.
Since $\hn(S^2)=1$ by \eqref{eq:sphere}, the product inequality is
strict, as stated in \eqref{eq:strict}.
\end{proof}

\section{The two corrections and the filtered obstruction}\label{sec:mechanism}
The two exact terms in \eqref{eq:trace-family} have distinct roles.
The weight-three term makes the quadratic sphere correction admissible.
The weight-zero term removes an obstruction at the other end of the
syzygy chain. We first characterize the parameter family, then describe
the resulting obstruction in a sphere-weight filtration.

\subsection{Necessity of both parameters}
\begin{proposition}\label{prop:parameter-locus}
For $\alpha,\beta\in\Q$, the following are equivalent:
\begin{enumerate}
\item $\alpha\beta\ne0$;
\item the ideal $J_{\alpha,\beta}$ of \eqref{eq:final-ideal} is acyclic;
\item there is an acyclic graded differential ideal $I'$ with
$\m_{\A}^{12}\subset I'\subset O_{\eps_{\alpha,\beta}}$.
\end{enumerate}
\end{proposition}
\begin{proof}
Theorem~\ref{thm:product-acyclicity} proves (1)$\Rightarrow$(2).
The implication (2)$\Rightarrow$(3) follows because
$\A\ctr\eps_{\alpha,\beta}\subset\E_{\alpha,\beta}$.
For the converse, we exhibit an obstruction on each coordinate axis.

Suppose $\beta=0$. The cycle $x^{11}u$ has length twelve. If $y\in\A^{23}$ satisfies $Dy=x^{11}u$, the weight-one component of $y$ is $u\pi$
with $D\pi=x^{11}$ and $\pi\in\N^{21}$.
Formula~\eqref{eq:all-primitives} gives
\[
 \pi=x^5b_0+\mu_0x^2\sigma_0+aP_a,
 \qquad P_a\in\Span_{\Q}\{x^9,x^6z,x^3z^2,z^3\}.
\]
As in \eqref{eq:uncorrected-obstruction}, the weight-one part of
$\eps_{\alpha,0}$ gives $\lambda(\pi az^2b_3)=3$.
The weight-zero part of $y$ contributes nothing: the only support of
$D^*\tau_0$ containing both $a$ and $b_3$ has a single factor $z$, so
it cannot contain the multiplier $az^2b_3$. All components of weight at
least two are invisible to $\eps_{\alpha,0}$. Hence
\[
 \eps_{\alpha,0}(yaz^2b_3)=3.
\]
No primitive belongs to $O_{\eps_{\alpha,0}}$, excluding (3).

Now suppose $\alpha=0$, with $\beta$ arbitrary. The length-twelve cycle
$c_{\mathrm{top}}=ax^5z^5u$ is exact in $\A$, since
\[
 D(-ax^5z^2b_3u)=c_{\mathrm{top}}.
\]
Set $\chi=[ax^5z^5u]^\vee$. Transposing the differential gives
\[
 D^*\chi=6[x^3z^6u]^\vee-[axz^4b_1u]^\vee
                     -[ax^3z^3b_2u]^\vee-[ax^5z^2b_3u]^\vee.
\]
The terms with $\beta$ in \eqref{eq:Psi} contain no $b_0$, and
\eqref{eq:Phi} has no such term either. Consequently
\begin{equation}\label{eq:upper-axis-obstruction}
 b_0\ctr\eps_{0,\beta}=-3D^*\chi.
\end{equation}
If a primitive $y$ of $c_{\mathrm{top}}$ belonged to
$O_{\eps_{0,\beta}}$, equation~\eqref{eq:upper-axis-obstruction} would
give
\[
 1=\chi(c_{\mathrm{top}})=\chi(Dy)
               =-\tfrac13\eps_{0,\beta}(yb_0)=0,
\]
a contradiction. Condition (3) fails on both axes, which completes the equivalence.
\end{proof}

Proposition~\ref{prop:parameter-locus} concerns orientation functionals
on the same minimal algebra. The homology nilpotency of $\A$ remains
eleven on the axes; a functional on an axis simply cannot contain an
optimal acyclic ideal in its radical. The criterion works on the whole
open set $\alpha\beta\ne0$, with all coefficients defined over $\Q$.

\subsection{The first obstruction dies at sphere weight two}
Fix $\eps=\eps_{1,1}$ and $J=J_{1,1}$. The first element of
\eqref{eq:five-elements} is
\begin{equation}\label{eq:quadratic-primitive}
 \rho_1=x^5b_0+3u^2\sigma_0=x^5b_0+D(3v\sigma_0),
 \qquad D\rho_1=x^{11}.
\end{equation}
It belongs not only to $O_\eps$ but also to $J$.
Indeed, it has one letter $c_i$ and no $a$, so three of the four
odd-letter types of $\F$ in \eqref{eq:final-ideal} vanish on it.
For the fourth, $B(\xi)(\rho_1)=\xi(x^{11})=0$ for every
$\xi\in\W_{10}$. Its membership in $O_\eps$ was proved in
Proposition~\ref{prop:five-elements}. Thus $\rho_1\in\E_{1,1}^\perp$.
In particular, $u\rho_1\in J$ is a primitive of the length-twelve cycle
$x^{11}u$.

Let $\A_w$ denote the sphere-weight-$w$ component of $\A$.
For either $\mathcal{L}=J$ or $\mathcal{L}=O_\eps$, define
\begin{equation}\label{eq:weight-filtration}
 F_p\mathcal{L}=\mathcal{L}\cap\bigoplus_{0\leq w\leq p}\A_w,
 \qquad F_{-1}\mathcal{L}=0.
\end{equation}
The subspaces in \eqref{eq:weight-filtration} are subcomplexes because
$D$ preserves the grading \eqref{eq:sphere-weight}. The
filtration is exhaustive and bounded in each cohomological degree.
Its intersections need not decompose as sums of weight-homogeneous
elements of $\mathcal{L}$; this is exactly where higher differentials can occur.

\begin{proposition}\label{prop:weight-two}
For both choices of $\mathcal{L}$,
\begin{equation}\label{eq:weight-persistence}
 [x^{11}]\ne0\text{ in }H^{22}(F_0\mathcal{L}),\qquad
 [x^{11}]\ne0\text{ in }H^{22}(F_1\mathcal{L}),\qquad
 [x^{11}]=0\text{ in }H^{22}(F_2\mathcal{L}).
\end{equation}
With the convention $E_0^{p,q}=F_p\mathcal{L}^{p+q}/F_{p-1}\mathcal{L}^{p+q}$,
the differentials have bidegree $(-r,r+1)$:
\[
 d_r:E_r^{p,q}\longrightarrow E_r^{p-r,q+r+1}.
\]
In particular, the spectral sequence has a nonzero differential
\begin{equation}\label{eq:weight-d2}
 d_2:E_2^{2,19}\longrightarrow E_2^{0,22},\qquad
 d_2[\operatorname{in}_2(\rho_1)]=[x^{11}].
\end{equation}
\end{proposition}
\begin{proof}
Equation~\eqref{eq:quadratic-primitive} shows that
$\rho_1\in F_2\mathcal{L}$ and $D\rho_1=x^{11}\in F_0\mathcal{L}$.
It proves the last assertion of \eqref{eq:weight-persistence}.
Suppose $y\in F_1\mathcal{L}^{21}$ were a primitive of $x^{11}$.
Since $v$ has weight two, the ambient decomposition has the form
$y=y_0+uy_1$, with $y_0\in\N^{21}$ and $y_1\in\N^{19}$.
Comparison of weight zero gives $Dy_0=x^{11}$.
The functional $\eps$ has weights zero, one, and three, but no
weight-two component. Therefore
\[
 0=\eps(yaz^2b_3u)=\lambda(y_0az^2b_3)=3,
\]
where the first equality uses $y\in \mathcal{L}\subset O_\eps$, and the last
uses \eqref{eq:all-primitives}. This contradiction proves nonvanishing
in $H^{22}(F_1\mathcal{L})$, and hence also in $H^{22}(F_0\mathcal{L})$.

The initial form of $\rho_1$ in filtration degree two is represented
in the ambient associated graded algebra by $3u^2\sigma_0$.
Since its differential has filtration degree zero, it survives the
first two associated-graded steps. It cannot be a boundary at either
of those steps: subtracting the corresponding actual boundary would
leave a primitive of $x^{11}$ in $F_1\mathcal{L}$, which has just been excluded.
Likewise, $[x^{11}]$ is nonzero on the second page, since a boundary
reaching it from filtration degree zero or one would give such a
primitive. The component of $D\rho_1$ two filtration levels below the
initial form is $x^{11}$, proving \eqref{eq:weight-d2}.
\end{proof}

An equivalent formulation uses the Rees complex
\[
 \mathcal R(\mathcal{L})=\bigoplus_{p\geq0}\hbar^pF_p\mathcal{L},
\]
where the parameter $\hbar$ has cohomological degree zero and $D$ acts on
coefficients. Equation~\eqref{eq:quadratic-primitive} gives
$D(\hbar^2\rho_1)=\hbar^2x^{11}$. By \eqref{eq:weight-persistence}, neither
$x^{11}$ nor $\hbar x^{11}$ is a boundary in $\mathcal R(\mathcal{L})$.
The submodule generated by $[x^{11}]$ in $H(\mathcal R(\mathcal{L}))$ is
therefore isomorphic to $\Q[\hbar]/(\hbar^2)$.

\subsection{The quotient obstruction and the Postnikov data}\label{sec:postnikov}
The product does not kill an ordinary $k$-invariant of $N$.
The projection onto each factor has a section, and the induced map on
rational cohomology is injective. In the pure description, the action
$\theta[z]=[x^2]\ne0$ in $H_0(\K_N)$ remains nonzero after tensoring
with the sphere model. The nontrivial twist in \eqref{eq:D-split} is
retained.

What changes is the ability to choose primitives inside one
multiplicatively closed ideal. The term $D(3v\sigma_0)$ changes a
primitive by an exact element of the ambient algebra while canceling
its forbidden orientation pairing. The last syzygy and the elements
$\rho_3,\rho_4,\rho_5$ supply the remaining compatibility relations.
The residue inclusion \eqref{eq:residue-target} turns these local
relations into a condition on all polynomial multipliers, and
\eqref{eq:comparison-diagram} turns that condition into acyclicity.
The filtered class of Proposition~\ref{prop:weight-two} is relative to
the specified trace and ideal; the recognition theorem is the separate
global statement that controls their entire quotient.

\section{Repeated products and further questions}\label{sec:consequences}
The first sphere factor removes the quotient defect of the manifold
$M$, while every further factor $S^2$ increases both invariants by one.

\begin{corollary}\label{cor:iterated}
For the manifold $M$ in Theorem~\ref{thm:main} and every integer $s\geq1$,
\[
 \hn\bigl(M\times(S^2)^s\bigr)
   =\cat_0\bigl(M\times(S^2)^s\bigr)=10+s.
\]
These are simply connected, closed smooth manifolds of dimension $53+2s$.
\end{corollary}
\begin{proof}
The product formula \eqref{eq:cat-product} gives rational category
$10+s$, which is a lower bound for homology nilpotency by
\eqref{eq:comparisons}. Apply \eqref{eq:hn-product} to
$(M\times S^2)\times(S^2)^{s-1}$ and use
Theorem~\ref{thm:main} and \eqref{eq:sphere}; the upper bound is
$11+(s-1)=10+s$. The dimension and manifold assertions follow from the
ordinary product construction.
\end{proof}

The example raises several questions about the role of the sphere
relation. For which dimensions $k$ can $\hn(X\times S^k)=\hn(X)$
hold for a noncontractible rationally elliptic space? Can the phenomenon
occur with fewer minimal generators or in smaller formal dimension?
The adjacent syzygies of $(t,z)^3$ enter the present construction;
other powers $(t,z)^m$ provide a natural class in which to study the
residue condition.

The recognition theorem assumes pure Poincaré duality with orientation
of Koszul degree two. It is natural to ask for an analogue in higher
Koszul degree, where contractions of the orientation occupy more than
three pure degrees. The central issue is whether polynomial residues
can still control all the ambient boundaries needed to reconstruct the
extension by a closed odd generator.

For general products, is
$\hn(X\times Y)\geq\max\{\hn(X),\hn(Y)\}$ always true?
How large can the \emph{additivity deficit}
$\hn(X)+\hn(Y)-\hn(X\times Y)$ be among rationally elliptic spaces,
or among simply connected closed smooth manifolds?
Theorem~\ref{thm:main} gives deficit one with a single sphere factor;
the possible size of the deficit and its relation to the quotient
defects of the factors remain distinct questions.

\makeatletter
\@setaddresses
\global\let\addresses\@empty
\makeatother

\begingroup
\makeatletter
\let\@tocwrite\@gobbletwo
\makeatother
\section*{Declaration of generative AI and AI-assisted technologies
in the manuscript preparation process}
During the preparation of this manuscript, the author used ChatGPT (OpenAI)
to assist with literature searches, language editing, \LaTeX\ preparation,
and review of mathematical arguments and exact
coefficient checks. The author takes full responsibility for all statements, proofs, references,
and conclusions, and for the article as a whole.
\par
\endgroup

\par\vfill

\begin{thebibliography}{99}

\bibitem{Barge}
J.~Barge,
\emph{Structures différentiables sur les types d'homotopie rationnelle simplement connexes},
Ann. Sci. Éc. Norm. Supér. (4) \textbf{9} (1976), no.~4, 469--501.
\doi{10.24033/asens.1315}.

\bibitem{Carrasquel}
J.~G.~Carrasquel-Vera,
\emph{Computations in rational sectional category},
Bull. Belg. Math. Soc. Simon Stevin \textbf{22} (2015), no.~3, 455--469.
\doi{10.36045/bbms/1442364592}.

\bibitem{CorneaCone}
O.~Cornea,
\emph{Cone-length and Lusternik--Schnirelmann category},
Topology \textbf{33} (1994), no.~1, 95--111.
\doi{10.1016/0040-9383(94)90037-X}.

\bibitem{Cornea}
O.~Cornea,
\emph{There is just one rational cone-length},
Trans. Amer. Math. Soc. \textbf{344} (1994), no.~2, 835--848.
\doi{10.1090/S0002-9947-1994-1260200-X}.

\bibitem{CFL}
O.~Cornea, Y.~Félix, and J.-M.~Lemaire,
\emph{Rational category and cone length of Poincaré complexes},
Topology \textbf{37} (1998), no.~4, 743--748.
\doi{10.1016/S0040-9383(97)00060-8}.

\bibitem{CLOT}
O.~Cornea, G.~Lupton, J.~Oprea, and D.~Tanré,
\emph{Lusternik--Schnirelmann Category}, Mathematical Surveys and Monographs,
vol.~103, American Mathematical Society, Providence, RI, 2003.
\doi{10.1090/surv/103}.

\bibitem{FH}
Y.~Félix and S.~Halperin,
\emph{Rational L.-S. category and its applications},
Trans. Amer. Math. Soc. \textbf{273} (1982), no.~1, 1--37.
\doi{10.1090/S0002-9947-1982-0664027-0}.

\bibitem{FHL}
Y.~Félix, S.~Halperin, and J.-M.~Lemaire,
\emph{The rational LS category of products and of Poincaré duality complexes},
Topology \textbf{37} (1998), no.~4, 749--756.
\doi{10.1016/S0040-9383(97)00061-X}.

\bibitem{FHT}
Y.~Félix, S.~Halperin, and J.-C.~Thomas,
\emph{Rational Homotopy Theory}, Graduate Texts in Mathematics, vol.~205,
Springer, New York, 2001.
\doi{10.1007/978-1-4613-0105-9}.

\bibitem{Ganea}
T.~Ganea,
\emph{Some problems on numerical homotopy invariants},
in \emph{Symposium on Algebraic Topology} (Battelle Seattle Res. Center,
Seattle, Wash., 1971), Lecture Notes in Mathematics, vol.~249,
Springer, Berlin, 1971, pp.~23--30.
\doi{10.1007/BFb0060892}.

\bibitem{HL}
S.~Halperin and J.-M.~Lemaire,
\emph{Notions of category in differential algebra},
in \emph{Algebraic Topology---Rational Homotopy} (Y.~Félix, ed.),
Lecture Notes in Mathematics, vol.~1318, Springer, Berlin, 1988,
pp.~138--154. \doi{10.1007/BFb0077800}.

\bibitem{Hess}
K.~P.~Hess,
\emph{A proof of Ganea's conjecture for rational spaces},
Topology \textbf{30} (1991), no.~2, 205--214.
\doi{10.1016/0040-9383(91)90006-P}.

\bibitem{Iwase}
N.~Iwase,
\emph{Ganea's conjecture on Lusternik--Schnirelmann category},
Bull. London Math. Soc. \textbf{30} (1998), no.~6, 623--634.
\doi{10.1112/S0024609398004548}.

\bibitem{Jessup}
B.~Jessup,
\emph{Rational L-S category and a conjecture of Ganea},
J. Pure Appl. Algebra \textbf{65} (1990), no.~1, 57--67.
\doi{10.1016/0022-4049(90)90100-V}.

\bibitem{LS}
P.~Lambrechts and D.~Stanley,
\emph{Poincaré duality and commutative differential graded algebras},
Ann. Sci. Éc. Norm. Supér. (4) \textbf{41} (2008), no.~4, 497--511.
\doi{10.24033/asens.2074}.

\bibitem{Parent}
P.-E.~Parent,
\emph{Homology nilpotency in rational homotopy theory: Cell attachments,
retractive towers, divergence, and stabilization},
preprint, September~2, 2026.
\href{https://arxiv.org/abs/2609.02980v1}{arXiv:2609.02980v1 [math.AT]}.

\bibitem{Sullivan}
D.~Sullivan,
\emph{Infinitesimal computations in topology},
Publ. Math. IHÉS \textbf{47} (1977), 269--331.
\doi{10.1007/BF02684341}.

\bibitem{Toomer}
G.~H.~Toomer,
\emph{Lusternik-Schnirelmann category and the Moore spectral sequence},
Math. Z. \textbf{138} (1974), no.~2, 123--143.
\doi{10.1007/BF01214229}.

\end{thebibliography}
\end{document}